\documentclass[a4paper,11pt]{article}

\usepackage{amsmath}
\usepackage{amsthm}
\usepackage{amssymb}
\usepackage{mathrsfs}
\usepackage{graphicx}
\usepackage{hyperref}

\usepackage{abstract} 

\newtheorem{thm}{Theorem}[section]
\newtheorem{cor}[thm]{Corollary}

\newtheorem{lemma}[thm]{Lemma}
\newtheorem{prop}[thm]{Proposition}

\theoremstyle{definition}
\newtheorem{definition}[thm]{Definition}
\newtheorem{ex}[thm]{Example}
\newtheorem{remark}[thm]{Remark}
\newtheorem{question}[thm]{Question}
\newtheorem{conj}[thm]{Conjecture}

\title{Hyperplanes of Squier's cube complexes}
\date{\today}
\author{Anthony Genevois}

\begin{document}

\maketitle

\begin{abstract}
To any semigroup presentation $\mathcal{P}= \langle \Sigma \mid \mathcal{R} \rangle$ and base word $w \in \Sigma^+$ may be associated a nonpositively curved cube complex $S(\mathcal{P},w)$, called a \emph{Squier complex}, whose underlying graph consists of the words of $\Sigma^+$ equal to $w$ modulo $\mathcal{P}$ where two such words are linked by an edge when one can be transformed into the other by applying a relation of $\mathcal{R}$. A group is a \emph{diagram group} if it is the fundamental group of a Squier complex. In this paper, we describe hyperplanes in these cube complexes. As a first application, we determine exactly when $S(\mathcal{P},w)$ is a special cube complex, as defined by Haglund and Wise, so that the associated diagram group embeds into a right-angled Artin group. A particular feature of Squier complexes is that the intersections of hyperplanes are ``ordered" by a relation $\prec$. As a strong consequence on the geometry of $S(\mathcal{P},w)$, we deduce, in finite dimensions, that its univeral cover isometrically embedds into a product of finitely-many trees with respect to the combinatorial metrics; in particular, we notice that (often) this allows to embed quasi-isometrically the associated diagram group into a product of finitely-many trees, giving information on its asymptotic dimension and its uniform Hilbert space compression. Finally, we exhibit a class of hyperplanes inducing a decomposition of $S(\mathcal{P},w)$ as a graph of spaces, and a fortiori a decomposition of the associated diagram group as a graph of groups, giving a new method to compute presentations of diagram groups. As an application, we associate a semigroup presentation $\mathcal{P}(\Gamma)$ to any finite \emph{interval graph} $\Gamma$, and we prove that the diagram group associated to $\mathcal{P}(\Gamma)$ (for a given base word) is isomorphic to the right-angled Artin group $A(\overline{\Gamma})$. This result has many consequences on the study of subgroups of diagram groups. In particular, we deduce that, for all $n \geq 1$, the right-angled Artin group $A(C_n)$ embeds into a diagram group, answering a question of Guba and Sapir.
\end{abstract}

\vspace{-0.5cm}

\tableofcontents

\section{Introduction}

\noindent
Given a class of groups, an interesting question raised by Wise is: may these groups be \textit{cubulated}, that is, do they act nicely on a CAT(0) cube complex? A positive answer yields interesting properties, depending on the action we find. In particular, cubulating hyperbolic 3-manifold groups was the key point in proving the virtual Haken conjecture \cite{MR3104553}. Known cubulated groups include, in particular, Coxeter groups \cite{MR1983376}, Artin groups of type FC \cite{MR1303028}, small cancellation groups \cite{MR2053602}, one-relator groups with torsion \cite{MR3118410} and free-by-cyclic groups \cite{MR3320891}. The so-called \textit{diagram groups}, mainly studied by Guba and Sapir, were cubulated by Farley in \cite{MR1978047}. 

\medskip \noindent
A simple definition is the following: Let $\mathcal{P}= \langle \Sigma \mid \mathcal{R} \rangle$ be a semigroup presentation and $w \in \Sigma^+$ a base word. We define the \textit{Squier complex} $S(\mathcal{P},w)$ as the cube complex whose vertices are the words of $\Sigma^+$ equal to $w$ modulo $\mathcal{P}$; whose edges, written $(a,u \to v,b)$, link two words $aub$ and $avb$ if $u=v \in \mathcal{R}$; whose $n$-cubes are similarly associated to the notation 
\begin{center}
$(a_1,u_1 \to v_1, \ldots, a_n,u_n \to v_n, a_{n+1})$.
\end{center}
Then, the \textit{diagram group} $D(\mathcal{P},w)$ is defined as the fundamental group of $S(\mathcal{P},w)$.

\medskip \noindent
Although Squier complexes turn out to be nonpositively curved, diagram groups have not been studied from the point of view of CAT(0) cube complexes. Motivated by our previous article \cite{arXiv:1505.02053} in which we proved that a diagram group not containing $\mathbb{Z}^2$ is free, answering a question of Guba and Sapir, we hope this approach may be fruitful. In this article, we study Squier complexes as nonpositively curved cube complexes, mainly by considering their hyperplanes. 

\medskip \noindent
The first question we are interested in is (see Section 2 for precise definitions):

\begin{question}
When is a diagram group the fundamental group of (compact) special cube complex?
\end{question}

\noindent
According to \cite{MR2377497}, consequences of this property include linearity and separability of some subgroups. Therefore, a natural problem is to determine when Squier complexes are special. A precise answer is given by our first result:

\begin{thm}\label{mainspecial}
Let $\mathcal{P}= \langle \Sigma \mid \mathcal{R} \rangle$ be a semigroup presentation and $w_0 \in \Sigma^+$ a base word. Then, the following assertions are equivalent:
\begin{itemize}
	\item[$(i)$] $S(\mathcal{P},w_0)$ is clean,
	\item[$(ii)$] $S(\mathcal{P},w_0)$ has no self-intersecting hyperplanes,
	\item[$(iii)$] there are no words $a,b,p \in \Sigma^+$ such that $w_0=ab$, $a=ap$ and $b=pb$ modulo $\mathcal{P}$ with $[p]_{\mathcal{P}} \neq \{ p\}$.
\end{itemize}
Moreover, $S(\mathcal{P},w_0)$ is special if and only if it satisfies the conditions above and the following one:
\begin{itemize}
	\item[$(iv)$] there are no words $a,u,v,w,b,p,q, \xi \in \Sigma^+$ such that $w_0=auvwb$, $au=au(v \xi)$, $wb=( \xi v)wb$ modulo $\mathcal{P}$ and $uv=p,vw=q \in \mathcal{R}$.
\end{itemize}
\end{thm}

\noindent
In particular, any compact Squier complex is special. Another geometric property specific to Squier complexes is that the intersections of hyperplanes are ``ordered". Roughly speaking, we say that, inside a square
$(a,u \to v, b, p \to q,c)$,
the hyperplane $J_1$ dual to the edge $(a, u \to v, bpc)$ meets the hyperplane $J_2$ dual to the edge $(aub,p \to q,c)$ \textit{by the left}; we note $J_1 \prec J_2$.  

\begin{prop}
The relation $\prec$ satisfies the following properties:
\begin{itemize}
	\item[$\bullet$] If $J_1 \prec J_2$ and $J_2 \prec J_3$, then $J_1 \prec J_3$.
	\item[$\bullet$] $J_1$ and $J_2$ are comparable with respect to $\prec$ if and only if they intersect. 
	\item[$\bullet$] $\mathrm{max} \{ n \geq 0 \mid \text{there exist} \ J_1, \ldots, J_n \ \text{such that} \ J_1 \prec \cdots \prec J_n \}= \dim S(\mathcal{P},w)$.
\end{itemize}
\end{prop}

\noindent
As a corollary, it is not difficult to deduce that the transversality graphs of Squier complexes do not contain induced cycles of odd length greater than five (Corollary \ref{cycle}), restricting the class of Squier complexes among nonpositively curved cube complexes. In the finite-dimensional case, we use the relation $\prec$ to prove:

\begin{thm}
Let $X(\mathcal{P},w)$ denote the universal cover of $S(\mathcal{P},w)$. If $d= \dim S(\mathcal{P},w)$ is finite, then $X(\mathcal{P},w)$ isometrically embeds into a product of $d$ trees with respect to the combinatorial metrics.
\end{thm}

\noindent
Thus, the second question we are interested in is:

\begin{question}
When does a finitely-generated diagram group quasi-isometrically embed into the product of finitely-many trees?
\end{question}

\noindent
A positive answer to this question gives information on the asymptotic geometry of the group: it bounds the asymptotic dimension, the dimensions of the asymptotic cones, and the uniform Hilbert space compression. Many groups are known to be quasi-isometrically embeddable into a product of finitely-many trees, such as hyperbolic groups \cite{MR2308852}, some relatively hyperbolic groups \cite{MR3073916}, mapping class groups \cite{arXiv:1207.2132} and virtually special groups. On the other hand, Thompson's group $F$, the discrete Heisenberg group and wreath products are known for not satisfying this property \cite{MR1883722}. In particular, since Thompson's group $F$ and several wreath products are diagram groups, the property we are considering does not hold for all diagram groups. We deduce a partial answer to the question above thanks to the previous theorem, provided that Property B, introduced in \cite{MR2271228}, is satisfied (although we believe that it is always the case):

\begin{thm}\label{mainQItrees}
Suppose that $S(\mathcal{P},w)$ is finite-dimensional and $D(\mathcal{P},w)$ is finitely-generated. If $D(\mathcal{P},w)$ satisfies Property B, then it quasi-isometrically embeds into a product of $\dim S(\mathcal{P},w)$ trees.
\end{thm}

\noindent
Finally, we introduce a family of hyperplanes in $S(\mathcal{P},w)$ which induces a decomposition of $S(\mathcal{P},w)$ as a graph of spaces; see Theorem \ref{graphofspaces} for a precise statement. In particular, this gives a decomposition of the associated diagram group $D(\mathcal{P},w)$ as a graph of groups. In fact, we already used a similar splitting in \cite{arXiv:1505.02053}.

\medskip \noindent
As an application, we will identify new diagram groups among the class of right-angled Artin groups. 

\begin{definition}
Let $\Gamma$ be a simplicial graph. Let $V(\Gamma)$ (resp. $E(\Gamma)$) denote its set of vertices (resp. edges). Then, we define the \textit{right-angled Artin group} $A(\Gamma)$ by the presentation
\begin{center}
$A(\Gamma)= \langle v \in V(\Gamma) \mid [u,v]=1, \ (u,v) \in E(\Gamma) \rangle$.
\end{center}
\end{definition}

\noindent
Determining for which graph $\Gamma$ the right-angled Artin group $A(\Gamma)$ is a diagram group is a wide open problem. However, Guba and Sapir proved the following results:

\begin{thm}\label{A(T)}
\emph{\cite[Theorem 7.8]{MR2193191}} Let $T$ be a finite tree. Then $A(T)$ is a diagram group.
\end{thm}

\begin{thm}\label{A(Cn)}
\emph{\cite[Theorem 30]{MR1725439}} Let $C_n$ be a cycle of length $n$. Then $A(C_n)$ is not a diagram group when $n \geq 5$ is odd.
\end{thm}

\noindent
However, it is not even known whether $A(C_n)$ is a diagram group when $n \geq 6$ is even. In fact, Guba and Sapir asked \cite[Problem 7]{MR1725439} whether $A(C_n)$ may be a subgroup of a diagram group; we answer the question below. 

\medskip \noindent
In our application, we will be interested in the following specific family of graphs:

\begin{definition}\label{intervalgraph}
To any collection $\mathcal{C}$ of intervals on the real line is associated a graph $\Gamma(\mathcal{C})$ whose set of vertices is $\mathcal{C}$ and whose edges link two intersecting intervals. We say that $\Gamma(\mathcal{C})$ is an \textit{interval graph}.
\end{definition}

\noindent
For example, the graph $P_n$, which denotes a path of length $n$, is an interval graph. The figure below gives a collection of intervals whose associated interval graph is $P_5$.
\begin{center}
\includegraphics[scale=0.42]{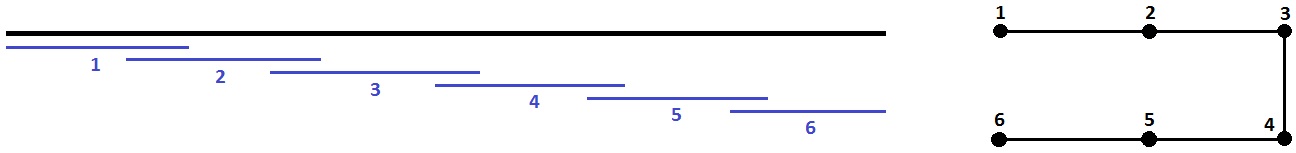}
\end{center}
In the following, we will be interested in \textit{complements} of interval graphs. Given a graph $\Gamma$, we define its \textit{complement} $\overline{\Gamma}$ as the graph whose vertices are the same as those of $\Gamma$ and whose edges link two vertices not linked by an edge in $\Gamma$. Therefore, the complement of $\Gamma(\mathcal{C})$ will be the graph whose set of vertices is $\mathcal{C}$ and whose edges link two disjoint intervals.

\begin{thm}\label{mainintervalgraphRAAG}
Let $\Gamma$ be a finite interval graph. Then the right-angled Artin group $A(\overline{\Gamma})$ is a diagram group.
\end{thm}

\noindent
It is not difficult to prove that the graph $\overline{\Gamma}$ cannot contain an induced path of length three, so that Theorem \ref{mainintervalgraphRAAG} and Theorem \ref{A(T)} essentially apply to different cases. As a first consequence, we answer \cite[Problem 7]{MR1725439}.

\begin{cor}\label{AC5}
For all $n \geq 1$, $A(C_n)$ is a subgroup of a diagram group.
\end{cor}

\noindent
\textbf{Proof.} Since $P_7$ is an interval graph, we deduce from Theorem~\ref{mainintervalgraphRAAG} that $A(\overline{P_7})$ is a diagram group. Then, according to \cite[Corollary 4.4]{MR3072113}, $A(C_n)$ embeds into $A(\overline{P_7})$ for all $n \geq 5$, so that the conclusion follows in this case. Finally, if $n \leq 4$, $A(C_n)$ is a diagram group. $\square$

\medskip \noindent
From the fact that $A(P_2(6))$ (see Example \ref{AP26}) embeds into a diagram group, it is deduced in \cite[Section 8]{MR2422070} that a diagram group may contain a hyperbolic surface group. This result is of interest because it proves that a diagram group may contain a subgroup whose first homology group has torsion or a hyperbolic subgroup which is not free; they are natural questions appearing in \cite{MR2193190}. Here, we are able to prove:

\begin{cor}\label{surfaces}
The fundamental group of a compact surface of even Euler characteristic at most $-2$ embeds into a diagram group. 
\end{cor}

\noindent
\textbf{Proof.} According to {\cite[Corollary 4.5]{MR3072113}, such a fundamental group embeds into $A(C_5)$. We conclude thanks to Corollary \ref{AC5}. $\square$

\medskip \noindent
In particular, every orientable surface group embeds into a diagram group, answering a question of Guba and Sapir \cite[$\S$ 17.3]{MR1396957}. Finally, we are able to give new counterexamples to the Subgroup Conjecture, ie., examples of subgroups of diagram groups which are not diagram groups themselves (see \cite{MR1725439} for the first known counterexamples).

\begin{cor}
For all odd $n \geq 5$, $A(C_n)$ embeds into a diagram group but is not a diagram group itself.
\end{cor}

\noindent
\textbf{Proof.} It is a consequence of Corollary~\ref{AC5} and Theorem~\ref{A(Cn)}. $\square$

\begin{cor}
The fundamental group of a hyperbolic closed surface of even Euler characteristic at most $-2$ embeds into a diagram group but is not a diagram group itself.
\end{cor}

\noindent
\textbf{Proof.} We already saw that such a fundamental group embeds into a diagram group. Moreover, it is a non-free hyperbolic group, so that it cannot be a diagram group according to \cite{arXiv:1505.02053}. $\square$

\medskip \noindent
The paper is organized as follows. In Section 2, we expose the preliminaries neeeded in the rest of the article; they concern diagram groups and cube complexes. In section 3, we describe hyperplanes in Squier complexes and we prove Theorem \ref{mainspecial}. In Section 4, we define the relation $\prec$, and we show how to use it to define the \textit{rank} of a hyperplane in the finite-dimensional case in order to finally prove Theorem \ref{mainQItrees}. In Section 5, we define left hyperplanes and exhibit a decomposition of the Squier complexes as graphs of spaces. As an application, we prove Theorem \ref{mainintervalgraphRAAG}. Finally, we conclude our article by some open questions in Section 6.

\section{Preliminaries}

\noindent
\textbf{2.1. Diagram groups.} We refer to \cite[$\S$3 and $\S$5]{MR1396957} for a detailed introduction to \textit{semigroup diagrams} and \textit{diagram groups}. 

For an alphabet $\Sigma$, let $\Sigma^+$ denote the free semigroup over $\Sigma$. If $\mathcal{P}= \langle \Sigma \mid \mathcal{R} \rangle$ is a semigroup presentation, where $\mathcal{R}$ is a set of pairs of words in $\Sigma^+$, the semigroup associated to $\mathcal{P}$ is the one given by the factor-semigroup $\Sigma^+ / \sim$ where $\sim$ is the smallest equivalent relation on $\Sigma^+$ containing $\mathcal{R}$. For convenience, we will assume that if $u=v \in \mathcal{R}$ then $v=u \notin \mathcal{R}$; in particular, $u=u \notin \mathcal{R}$.

A \textit{semigroup diagram over $\mathcal{P}$} is the analogue for semigroups to van Kampen diagrams for group presentations. Formally, it is a finite connected planar graph $\Delta$ whose edges are oriented and labelled by the alphabet $\Sigma$, satisfying the following properties:\\
\indent $\bullet$ $\Delta$ has exactly one vertex-source $\iota$ (which has no incoming edges) and exactly one vertex-sink $\tau$ (which has no outgoing edges);\\
\indent $\bullet$ the boundary of each cell has the form $pq^{-1}$ where $p=q$ or $q=p \in \mathcal{R}$;\\
\indent $\bullet$ every vertex belongs to a positive path connecting $\iota$ and $\tau$;\\
\indent $\bullet$ every positive path in $\Delta$ is simple.\\
In particular, $\Delta$ is bounded by two positive paths: the \textit{top path}, denoted $\text{top}(\Delta)$, and the \textit{bottom path}, denoted $\text{bot}(\Delta)$. By extension, we also define $\text{top}(\Gamma)$ and $\text{bot}(\Gamma)$ for every \textit{subdiagram} $\Gamma$. In the following, the notations $\text{top}(\cdot)$ and $\text{bot}(\cdot)$ will refer to the paths and to their labels. Also, a \textit{$(u,v)$-cell} (resp. a \textit{$(u,v)$-diagram}) will refer to a cell (resp. a semigroup diagram) whose top path is labelled by $u$ and whose bottom path is labelled by $v$.

Two words $w_1,w_2$ in $\Sigma^+$ are \textit{equal modulo $\mathcal{P}$} if their images in the semigroup associated to $\mathcal{P}$ coincide. In particular, there exists a \textit{derivation} from $w_1$ to $w_2$, i.e., a sequence of relations of $\mathcal{R}$ allowing us to transform $w_1$ into $w_2$. To any such derivation is associated a semigroup diagram, or more precisely a $(w_1,w_2)$-diagram, whose construction is clear from the example below. As in the case for groups, the words $w_1,w_2$ are equal modulo $\mathcal{P}$ if and only if there exists a $(w_1,w_2)$-diagram.

\begin{ex}
Let $\mathcal{P}= \langle a,b,c \mid ab=ba, ac=ca, bc=cb \rangle$ be a presentation of the free Abelian semigroup of rank three. In particular, the words $a^2bc$ and $caba$ are equal modulo $\mathcal{P}$, with the following possible derivation:
\begin{center}
$aabc \overset{(a,ab \to ba,c)}{\longrightarrow} abac \overset{(ab,ac \to ca, \emptyset)}{\longrightarrow} abca \overset{(a,bc \to cb,a)}{\longrightarrow} acba \overset{(\emptyset,ac \to ca,ba)}{\longrightarrow} caba$.
\end{center}
Then, the associated $(a^2bc,caba)$-diagram $\Delta$ is:
\begin{center}
\includegraphics[scale=0.6]{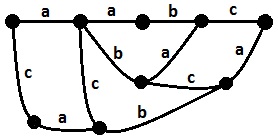}
\end{center}
On such a graph, the edges are oriented from left to right. Here, the diagram $\Delta$ has nine vertices, twelve edges and four cells; notice that the number of cells of a diagram corresponds to the length of the associated derivation. The paths $\text{top}(\Delta)$ and $\text{bot}(\Delta)$ are labelled respectively by $a^2bc$ and $caba$. 
\end{ex}

Since we are only interested in the combinatorics of semigroup diagrams, we will not distinguish isotopic diagrams. For example, the two diagrams below will be considered as equal.
\begin{center}
\includegraphics[scale=0.6]{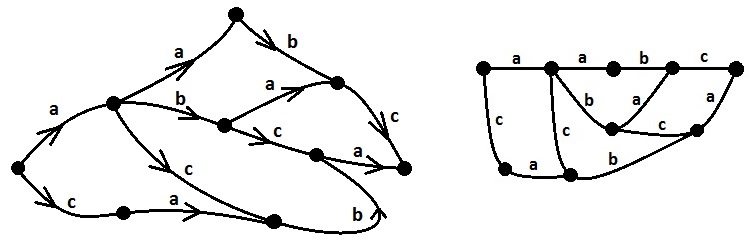}
\end{center}
If $w \in \Sigma^+$, we define the \textit{trivial diagram} $\epsilon(w)$ as the semigroup diagram without cells whose top and bottom paths, labelled by $w$, coincide. Any diagram without cells is trivial. A diagram with exactly one cell is \textit{atomic}.

If $\Delta_1$ is a $(w_1,w_2)$-diagram and $\Delta_2$ a $(w_2,w_3)$-diagram, we define the \textit{concatenation} $\Delta_1 \circ \Delta_2$ as the semigroup diagram obtained by identifying the bottom path of $\Delta_1$ with the top path of $\Delta_2$. In particular, $\Delta_1 \circ \Delta_2$ is a $(w_1,w_3)$-diagram. Thus, $\circ$ defines a partial operation on the set of semigroup diagrams over $\mathcal{P}$. However, restricted to the subset of $(w,w)$-diagrams for some $w \in \Sigma^+$, it defines a semigroup operation; such diagrams are called \textit{spherical with base $w$}. We also define the \textit{sum} $\Delta_1+ \Delta_2$ of two diagrams $\Delta_1,\Delta_2$ as the diagram obtained by identifying the rightmost vertex of $\Delta_1$ with the leftmost vertex of $\Delta_2$.
\begin{center}
\includegraphics[scale=0.6]{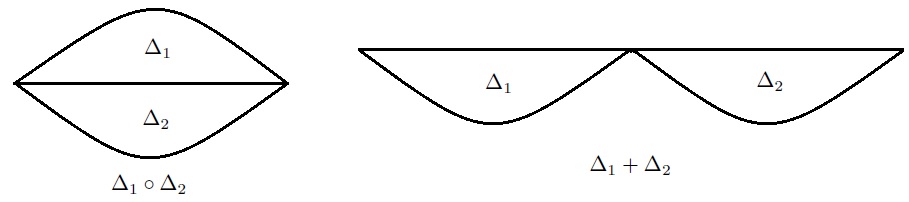}
\end{center}
Notice that any semigroup diagram can be viewed as a concatenation of atomic diagrams. In the following, if $\Delta_1,\Delta_2$ are two diagrams, we will say that $\Delta_1$ is a \textit{prefix} (resp. a \textit{suffix}) of $\Delta_2$ if there exists a diagram $\Delta_3$ satisfying $\Delta_2= \Delta_1 \circ \Delta_3$ (resp. $\Delta_2= \Delta_3 \circ \Delta_1$). Throughout this paper, the fact that $\Delta$ is a prefix of $\Gamma$ will be denoted by $\Delta \leq \Gamma$.

Suppose that a diagram $\Delta$ contains a $(u,v)$-cell and a $(v,u)$-cell such that the top path of the first cell is the bottom path of the second cell. Then, we say that these two cells form a \textit{dipole}. In this case, we can remove these two cells by first removing their common path, and then identifying the top path of the first cell with the bottom path of the second cell; thus, we \textit{reduce the dipole}. A diagram is called \textit{reduced} if it does not contain dipoles. By reducing dipoles, a diagram can be transformed into a reduced diagram, and a result of Kilibarda \cite[Theorem 2.1]{MR1448329} proves that this reduced form is unique. If $\Delta_1,\Delta_2$ are two diagrams for which $\Delta_1 \circ \Delta_2$ is well defined, let us denote by $\Delta_1 \cdot \Delta_2$ the reduced form of $\Delta_1 \circ \Delta_2$.

If $w \in \Sigma^+$, we define the \textit{diagram group} $D(\mathcal{P},w)$ as the set of reduced $(w,w)$-diagrams endowed with the product $\cdot$ we defined above. If $\Delta$ is a $(w_1,w_2)$-diagram, let $\Delta^{-1}$ denote the $(w_2,w_1)$-diagram obtained from $\Delta$ by a mirror reflection with respect to $\mathrm{top}(\Delta)$. It can be noticed that, if $\Delta$ is a spherical diagram, then $\Delta^{-1}$ is the inverse of $\Delta$ in the associated diagram group.

Although this definition of $D(\mathcal{P},w)$ does not seem to give much information on its group structure, it allows to define a class of \textit{canonical subgroups}. If $\Gamma$ is a $(w,u)$-diagram and if we write $u=x_1u_1 \cdots x_nu_nx_{n+1}$, where the $x_i$ and $u_i$ are (possibly empty) subwords of $u$, then the map
\begin{center}
$(U_1,\ldots, U_n) \mapsto \Gamma \cdot \left( \epsilon(x_1) + U_1 + \cdots + \epsilon(x_n)+ U_n+ \epsilon(x_{n+1}) \right) \cdot \Gamma^{-1}$
\end{center}
defines an embedding from $D(\mathcal{P},u_1) \times \cdots \times D(\mathcal{P},u_n)$ into $D(\mathcal{P},w)$.

\medskip \noindent
\textbf{2.2. Special cube complexes.} A \textit{cube complex} is a CW-complex constructed by gluing together cubes of arbitrary (finite) dimension by isometries along their faces. Furthermore, it is \textit{nonpositively curved} if the link of any of its vertices is a simplicial \textit{flag} complex (ie., $n+1$ vertices span a $n$-simplex if and only if they are pairwise adjacent), and \textit{CAT(0)} if it is nonpositively curved and simply-connected. See \cite[page 111]{MR1744486} for more information.

\medskip \noindent
A fundamental feature of cube complexes is the notion of \textit{hyperplane}. Let $X$ be a nonpositively curved cube complex. Formally, a \textit{hyperplane} $J$ is an equivalence class of edges, where two edges $e$ and $f$ are equivalent whenever there exists a sequence of edges $e=e_0,e_1,\ldots, e_{n-1},e_n=f$ where $e_i$ and $e_{i+1}$ are parallel sides of some square in $X$. Notice that a hyperplane is uniquely determined by one of its edges, so if $e \in J$ we say that $J$ is the \textit{hyperplane dual to $e$}. Geometrically, a hyperplane $J$ is rather thought of as the union of the \textit{midcubes} transverse to the edges belonging to $J$.
\begin{center}
\includegraphics[scale=0.4]{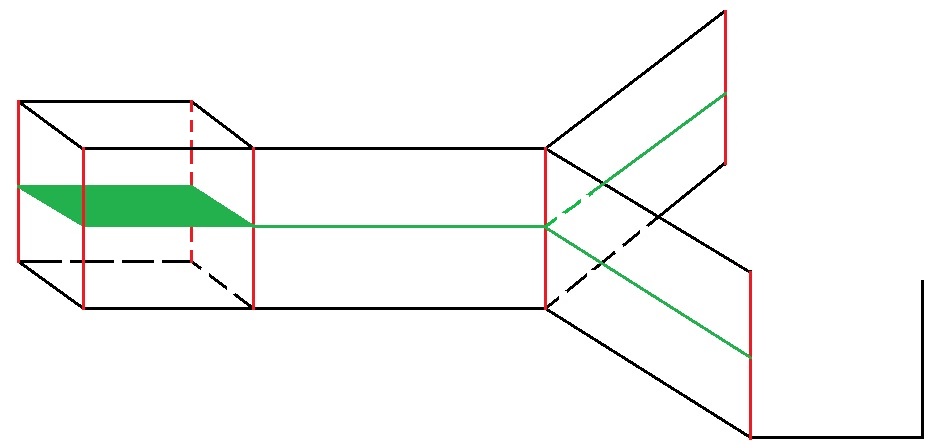}
\end{center}
Similarly, one may define \textit{oriented hyperplanes} as classes of oriented edges. If $J$ is the hyperplane dual to an edge $e$ and if we fix an orientation $\vec{e}$, we will note $\vec{J}$ the oriented hyperplane dual to $\vec{e}$. It may be thought of as an \textit{orientation} of $J$, and we will note $- \vec{J}$ the opposite orientation of $J$. 

\begin{definition}
Let $J$ be a hyperplane with a fixed orientation $\vec{J}$. We say that $J$ is:
\begin{itemize}
	\item \textit{2-sided} if $\vec{J} \neq - \vec{J}$,
	\item \textit{self-intersecting} if there exist two edges dual to $J$ which are non-parallel sides of some square,
	\item \textit{self-osculating} if there exist two oriented edges dual to $\vec{J}$ with the same initial points or the same terminal points, but which do not belong to a same square.
\end{itemize}
Moreover, if $H$ is another hyperplane, then $J$ and $H$ are:
\begin{itemize}
	\item \textit{transverse} if there exist two edges dual to $J$ and $H$ respectively which are non-parallel sides of some square,
	\item \textit{inter-osculating} if they are transverse, and if there exist two edges dual to $J$ and $H$ respectively with one common endpoint, but which do not belong to a same square.
\end{itemize}
\end{definition}
\noindent
Sometimes, one refers 1-sided, self-intersecting, self-osculating and inter-osculating hyperplanes as \textit{pathological configurations of hyperplanes}. The last three configurations are illustrated below.
\begin{figure}[h!]
\includegraphics[scale=0.21]{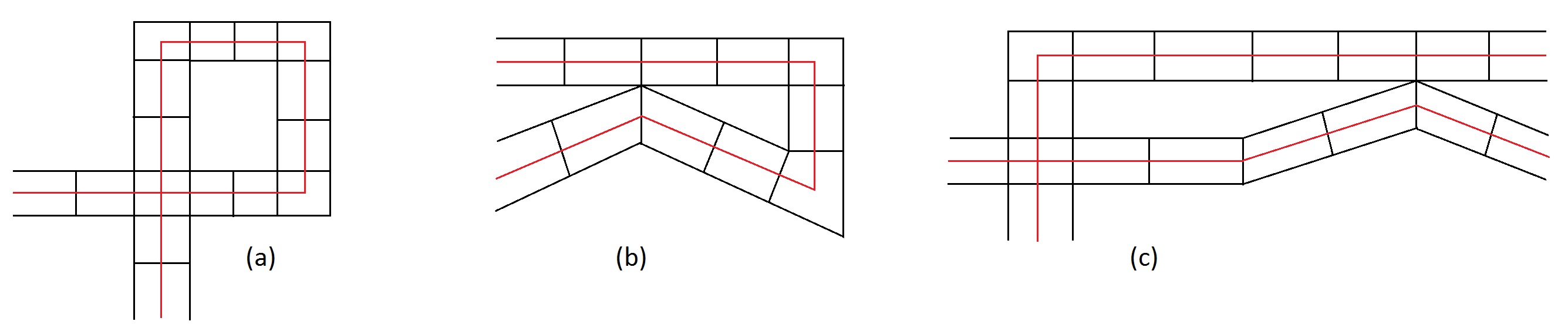}
\caption{(a) Self-intersection, (b) Self-osculation, (c) Inter-osculation.}
\end{figure}

\begin{definition}
A hyperplane is \textit{clean} if it is 2-sided and is neither self-intersecting nor self-osculating. A nonpositively curved cube complex is \textit{special} if its hyperplanes are clean and if it does not contain inter-osculating hyperplanes.
\end{definition}

Therefore, roughly speaking, special cube complexes are the cube complexes in which hyperplanes are ``well-behaved". They include CAT(0) cube complexes, in which hyperplanes satisfy the following properties:

\begin{thm}
Let $X$ be a CAT(0) cube complex. Then 
\begin{itemize}
	\item \emph{\cite[Example 3.3]{MR2377497}} $X$ is a special cube complex,
	\item \emph{\cite[Theorem 4.10]{MR1347406}} every hyperplane separates $X$ into two pieces, called \emph{halfspaces}.
\end{itemize}
\end{thm}

An important property of special cube complexes is that their fundamental groups embed into a right-angled Artin group. More precisely, we first define the graph we are interested in:

\begin{definition}
The \textit{transversality graph} of a cube complex $X$ is defined as the graph whose vertices are the hyperplanes of $X$ and whose edges link two transverse hyperplanes.
\end{definition}

Let $X$ be a special cube complex and let $\Gamma$ denote its transversality graph. Then, we can define a combinatorial map from $X$ to the \textit{Salvetti complex} $S(\Gamma)$ of $A(\Gamma)$ (roughly speaking, it is the CW-complex associated to the canonical presentation of $A(\Gamma)$ with additional cubes of dimensions $\geq 3$; it is a classifying space of $A(\Gamma)$) as follows: First, we fix an orientation of the hyperplanes of $X$; because these hyperplanes are 2-sided, it induces a well-defined orientation on the edges of $X$. Then, to any edge $e$ of $X$ is associated an oriented hyperplane, which is also a vertex of $\Gamma$, a generator of $A(\Gamma)$, and so naturally an (oriented) edge of the Salvetti complex $S(\Gamma)$; whence a map $\Psi$ from the edges of $X$ to the edges of $S(\Gamma)$. In fact, $\Psi$ may be extended into a combinatorial map $\Psi : X \to S(\Gamma)$ so that

\begin{thm}\label{localisometry}
\emph{\cite[Lemma 4.1]{MR2377497}} The cube complex $X$ is special if and only if $\Psi$ is a local isometry.
\end{thm}

\noindent
Because local isometries between nonpositively curved cube complexes are $\pi_1$-injective, we deduce

\begin{cor}\label{specialRAAG}
\emph{\cite[Theorem 4.4]{MR2377497}}The fundamental group of special cube complex embeds into a right-angled Artin groups.
\end{cor}

\noindent
In particular, we deduce that such a fundamental group is necessarily residually finite. In fact, when the cube complex is compact, it is possible to say more about separability properties.

\begin{definition}
Let $G$ be a group acting on a CAT(0) cube complex. A subgroup $H \leq G$ is convex-cocompact if there exists an $H$-invariant convex subcomplex $Y \subset X$ such that the action $H \curvearrowright Y$ is cocompact.
\end{definition}

\begin{thm}\label{convex-cocompact}
\emph{\cite[Corollary 7.9]{MR2377497}} Any convex-cocompact subgroup of the fundamental group of a compact special cube complex is separable.
\end{thm}

\noindent
Recall that a subgroup $H \leq G$ is \textit{separable} provided that, for all $g \in G \backslash H$, there exists a finite-index subgroup $K \leq G$ containing $H$ but not containing $g$.

\medskip \noindent
\textbf{2.3. Cubulation.} We now describe a method introduced by Sageev \cite{MR1347406, SageevCAT(0)}, called \textit{cubulation}, to construct CAT(0) cube complexes, which will be useful in Section 4. 

\begin{definition}
A \textit{pocset} $(\Sigma,<,^*)$ is a partially ordered set $(\Sigma,<)$ endowed with an involution $^*$ satisfying:
\begin{itemize}
	\item for all $A \in \Sigma$, $A$ and $A^*$ are not comparable (in particular, $A^* \neq A$),
	\item for all $A,B \in \Sigma$, $A< B$ if and only if $B^* < A^*$.
\end{itemize}
\end{definition}

\begin{definition}
Let $(\Sigma,<,^*)$ be a pocset. An \textit{ultrafilter} $\alpha$ is a set of subsets of $\Sigma$ satisfying:
\begin{itemize}
	\item for all $A \in \Sigma$, $\alpha$ contains exactly one element of $\{A,A^*\}$,
	\item for all $A,B \in  \Sigma$, if $B<A$ and $B \in \alpha$ then $A \in \alpha$.
\end{itemize}
Furthermore, $\alpha$ is a \textit{DCC ultrafilter} if every infinite descending chain in $\alpha$ is eventually constant.
\end{definition}

\noindent
Let $(\Sigma,<,^*)$ be a pocset. We define a cube complex $X(\Sigma)$ as follows:
\begin{itemize}
	\item the vertices are the DCC ultrafilters,
	\item two ultrafilters are linked by an edge if they differ by two subsets of $\Sigma$,
	\item we add $n$-cubes inductively as soon as possible, ie., we add 3-cubes as soon as the boundary of a 3-cube appear in the 2-skeleton, then we add 4-cubes as soon as the boundary of a 4-cube appear in the 3-skeleton, and so on.
\end{itemize}
We say that $X(\Sigma)$ is the cube complex constructed by \textit{cubulating} $\Sigma$.

\begin{thm}
Every connected component of $X(\Sigma)$ is a CAT(0) cube complex.
\end{thm}

\noindent
All the examples of pocsets we will consider come from the following one:

\begin{ex}
Let $X$ be a CAT(0) cube complex and $\mathcal{H}$ a collection of hyperplanes of $X$. The set $\Sigma$ of halfspaces delimited by the hyperplanes of $\mathcal{H}$ defines a pocset with respect to the inclusion $\subset$ and the complementary operation $^c$. Then, an ultrafilter $\alpha$ may be thought of as the choice of a halfspace for each hyperplane in $\mathcal{H}$. In particular, if $v \in X$ is a vertex, then the set of halfspaces of $\Sigma$ containing $v$ defines an ultrafilter, called \textit{principal}. Note that principal ultrafilters are DCC.
\end{ex}

In general, there is no canonical choice of a connected component of the cube complex constructed by cubulation. However, in the context of the previous example, we usually choose the connected component containing the principal ultrafilters (it can be shown that they all belong to the same component).

\medskip
Although a CAT(0) cube complex $X$ can be endowed with a CAT(0) metric, it is often more convenient to introduce a more ``combinatorial" distance. We define the \textit{combinatorial distance} $d_c$ on the set of vertices $X^{(0)}$ of $X$ as the graph metric associated to the 1-skeleton $X^{(1)}$. In fact, the combinatorial metric and the hyperplanes are linked together: it can be proved that the combinatorial distance between two vertices corresponds exactly to the number of hyperplanes separating them \cite[Theorem 2.7]{MR2413337}. This point of view allows us to link $X$ to the complexes constructed by cubulation from a collection of hyperplanes.

\begin{prop}\label{distancequotient}
Let $X$ be a CAT(0) cube complex, $\mathcal{H}$ a collection of hyperplanes, and let $X(\mathcal{H})$ denote the cube complex constructed by cubulation with respect to the pocset of halfspaces delimited by the hyperplanes of $\mathcal{H}$. Finally, let
\begin{center}
$\varphi : X^{(0)} \to X(\mathcal{H})^{(0)}$
\end{center}
be the natural map sending a vertex of $X$ to the principal ultrafilter it defines. Then, for any vertices $x,y \in X$, the combinatorial distance in $X(\mathcal{H})$ between $\varphi(x)$ and $\varphi(y)$ corresponds to the number of hyperplanes of $\mathcal{H}$ separating $x$ and $y$ in $X$.
\end{prop}

\medskip \noindent
\textbf{2.4. Squier and Farley complexes.} Let $\mathcal{P}= \langle \Sigma \mid \mathcal{R} \rangle$ be a semigroup presentation and $w \in \Sigma^+$ a base word. One way of obtaining information about the diagram group $D(\mathcal{P},w)$ is to describe it as the fundamental group of a cube complex.

More precisely, we define the \textit{Squier complex} $S(\mathcal{P})$ as the cube complex whose vertices are the words in $\Sigma^+$; whose (oriented) edges can be written as $(a,u \to v,b)$, where $u=v$ or $v=u$ belongs to $\mathcal{R}$, linking the vertices $aub$ and $avb$; and whose $n$-cubes similarly can be written as $(a_1,u_1 \to v_1, \ldots, a_n, u_n \to v_n, a_{n+1})$, spanned by the set of vertices $\{ a_1w_1 \cdots a_nw_na_{n+1} \mid w_i=v_i \ \text{or} \ u_i \}$.

Then, there is a natural morphism from the fundamental group of $S(\mathcal{P})$ based at $w$ to the diagram group $D(\mathcal{P},w)$. Indeed, a loop in $S(\mathcal{P})$ based at $w$ is just a series of relations of $\mathcal{R}$ applied to the word $w$ so that the final word is again $w$, and such a sequence may be encoded into a semigroup diagram. The figure below shows an example, where the semigroup presentation is $\mathcal{P}= \langle a,b,c \mid ab=ba, bc=cb,ac=ca \rangle$: 
\begin{center}
\includegraphics[scale=0.5]{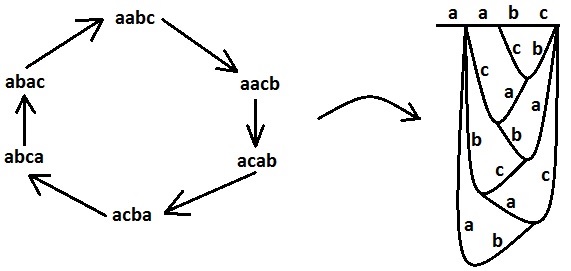}
\end{center}
Thus, this defines a map from the set of loops of $S(\mathcal{P})$ based at $w$ to the set of spherical semigroup diagrams. In fact, the map extends to a morphism which turns out to be an isomorphism:

\begin{thm}\label{iso}
\emph{\cite[Theorem 6.1]{MR1396957}} $D(\mathcal{P},w) \simeq \pi_1( S(\mathcal{P}),w)$.
\end{thm}

For convenience, $S(\mathcal{P},w)$ will denote the connected component of $S(\mathcal{P})$ containing $w$. Notice that two words $w_1,w_2 \in \Sigma^+$ are equal modulo $\mathcal{P}$ if and only if they belong to the same connected component of $S(\mathcal{P})$. Therefore, a consequence of Theorem \ref{iso} is:

\begin{cor}\label{morphism1}
If $w_1,w_2 \in \Sigma^+$ are equal modulo $\mathcal{P}$, then there exists a $(w_2,w_1)$-diagram $\Gamma$ and the map
\begin{center}
$\Delta \mapsto \Gamma \cdot \Delta \cdot \Gamma^{-1}$
\end{center}
induces an isomorphism from $D(\mathcal{P},w_1)$ to $D(\mathcal{P},w_2)$.
\end{cor}

\noindent
Another morphism between diagram groups which will be useful in Section 5 is:

\begin{lemma}\label{morphism2}
Let $u, v \in \Sigma^+$ be two words. Then the application
\begin{center}
$\Delta \mapsto \epsilon(u)+ \Delta$
\end{center}
induces a monomorphism from $D(\mathcal{P},v)$ into $D(\mathcal{P},uv)$.
\end{lemma}

\noindent
It can be proved that $S(\mathcal{P},w)$ is nonpositively curved so that its universal cover is CAT(0). In \cite{MR1978047}, Farley gives a construction of this cover.

A semigroup diagram is \textit{thin} whenever it can be written as a sum of atomic diagrams. We define the \textit{Farley complex} $X(\mathcal{P},w)$ as the cube complex whose vertices are the reduced semigroup diagrams $\Delta$ over $\mathcal{P}$ satisfying $\mathrm{top}(\Delta)=w$, and whose $n$-cubes are spanned by the vertices $\{ \Delta \cdot P \mid P \leq \Gamma \}$ for some vertex $\Delta$ and some thin diagram $\Gamma$ with $n$ cells. In particular, two diagrams $\Delta_1$ and $\Delta_2$ are linked by an edge if and only if there exists an atomic diagram $A$ such that $\Delta_1= \Delta_2 \cdot A$.

\begin{thm}
\emph{\cite[Theorem 3.13]{MR1978047}} $X(\mathcal{P},w)$ is a CAT(0) cube complex. Moreover it is complete, i.e., there is no increasing sequence of cubes in $X(\mathcal{P},w)$.
\end{thm}

\noindent
There is a natural action of $D(\mathcal{P},w)$ on $X(\mathcal{P},w)$, namely $(g, \Delta) \mapsto g \cdot \Delta$. Then

\begin{prop}\label{Farleyaction}
\emph{\cite[Theorem 3.13]{MR1978047}} The action $D(\mathcal{P},w) \curvearrowright X(\mathcal{P},w)$ is free. Moreover, it is properly discontinuous if $\mathcal{P}$ is a finite presentation, and it is cocompact if and only if the class $[w]_{\mathcal{P}}$ of words equal to $w$ modulo $\mathcal{P}$ is finite.
\end{prop}

\noindent
It is not difficult to prove that $\mathcal{P}$ may be supposed to be a finite presentation whenever $D(\mathcal{P},w)$ is finitely generated. Therefore, the action $D(\mathcal{P},w) \curvearrowright X(\mathcal{P},w)$ is often properly discontinuous.

\medskip \noindent
To conclude, we notice that the map $\Delta \mapsto \mathrm{bot}(\Delta)$ induces the universal covering $X(\mathcal{P},w) \to S(\mathcal{P},w)$ and that the action of $\pi_1(S(\mathcal{P},w))$ on $X(\mathcal{P},w)$ coincides with the natural action of $D(\mathcal{P},w)$. Precisely:

\begin{lemma}
\emph{\cite[Lemma 1.3.5]{arXiv:1505.02053}} The map $\Delta \mapsto \mathrm{bot}(\Delta)$ induces a cellular isomorphism from the quotient $X(\mathcal{P},w)/D(\mathcal{P},w)$ to $S(\mathcal{P},w)$.
\end{lemma}

\section{Specialness}

\subsection{Hyperplanes in Squier complexes}

\noindent
In this section, we fix a semigroup presentation $\mathcal{P}= \langle \Sigma \mid \mathcal{R} \rangle$ and a base word $w \in \Sigma^+$. If $(a, u \to v,b)$ is an edge in the Squier complex $S(\mathcal{P},w)$, $[a,u \to v,b]$ will denote the hyperplane dual to it. Below, we show how to use this notation to completely characterize the hyperplanes in $S(\mathcal{P},w)$ and we determine when two of them are transverse.

\begin{lemma}\label{hyperplansSquier}
Let $(p, a \to b , q)$ and $(r, c \to d ,s)$ be two edges in $S(\mathcal{P})$. Then the (oriented) hyperplanes $[p,a \to b,q]$ and $[r,c \to d,s]$ coincide if and only if $a=c$, $b=d$ in $\Sigma^+$ and $p=r$, $q=s$ modulo $\mathcal{P}$.
\end{lemma}

\noindent
\textbf{Proof.} First, suppose that $[p,a \to b,q]=[r,c \to d,s]$. We show the desired equalities by induction on the length $\ell$ of a path of parallel edges between $(p,a \to b,q)$ and $(r,c \to d,s)$. If $\ell=0$, then $(p,a \to b,q)=(r,c \to d,s)$ and there is nothing to prove. Suppose $\ell=n+1$. By the induction hypothesis, the $n$-th edge of our path can be written as $(p',a \to b,q')$ for some words $p',q' \in \Sigma^+$ satisfying $p'=p$ and $q'=q$ modulo $\mathcal{P}$; furthermore, this edge is parallel to $(r,c \to d,s)$, ie., these edges belong to a square
\begin{center}
$(x,u \to v,y,a \to b,q')$ or $(p',a \to b, x,u \to v, y)$,
\end{center}
where $xuy=p'$ modulo $\Sigma^+$ in the first case, and $xuy=q'$ modulo $\Sigma^+$ in the second case. Consequently,
\begin{center}
$(r, c \to d,s)= (xvy,a \to b,q')$ or $(p',a \to b,xvy)$
\end{center}
Thus, modulo $\mathcal{P}$, $r=p'=p$ and $s=q'=q$, and in $\Sigma^+$, $c=a$ and $d=b$.

\medskip \noindent
Conversely, we show that the edges $(p,a \to b,q)$ and $(r, a \to b,s)$ are dual to the same hyperplane provided that $p=r$ and $q=s$ modulo $\mathcal{P}$. In this situation, there exist two sequences
\begin{center}
$p=x_1,x_2,\ldots, x_n=r$ and $q=y_1,y_2, \ldots, y_m=s$
\end{center}
where the relations $x_i=x_{i+1}$ and $y_i=y_{i+1}$ belong to $\mathcal{R}$. We deduce the configuration:
\begin{center}
\includegraphics[scale=0.55]{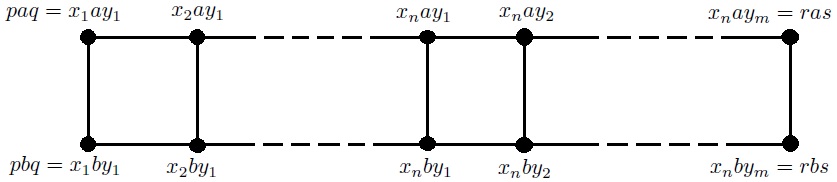}
\end{center}
Therefore, $[p,a \to b,q]=[r,a \to b,s]$. $\square$

\begin{lemma}\label{hyptransverseSquier}
Let $(a,u \to v,b)$ and $(c,p \to q,d)$ be two edges in $S(\mathcal{P})$. Then, the hyperplanes $[a, u \to v,b]$ and $[c, p \to q,d]$ are transverse if and only if there exists $y \in \Sigma^+$ satisfying $\left\{ \begin{array}{l} c=auy \\ b=ypd \end{array} \right.$ or $\left\{ \begin{array}{l} d=yub \\ a=cpy \end{array} \right.$ modulo $\mathcal{P}$. 
\end{lemma}

\noindent
\textbf{Proof.} Suppose that such a word $y$ exists. Then, we deduce one of the following two configurations, proving that the considered hyperplanes are transverse:
\begin{center}
\includegraphics[scale=0.4]{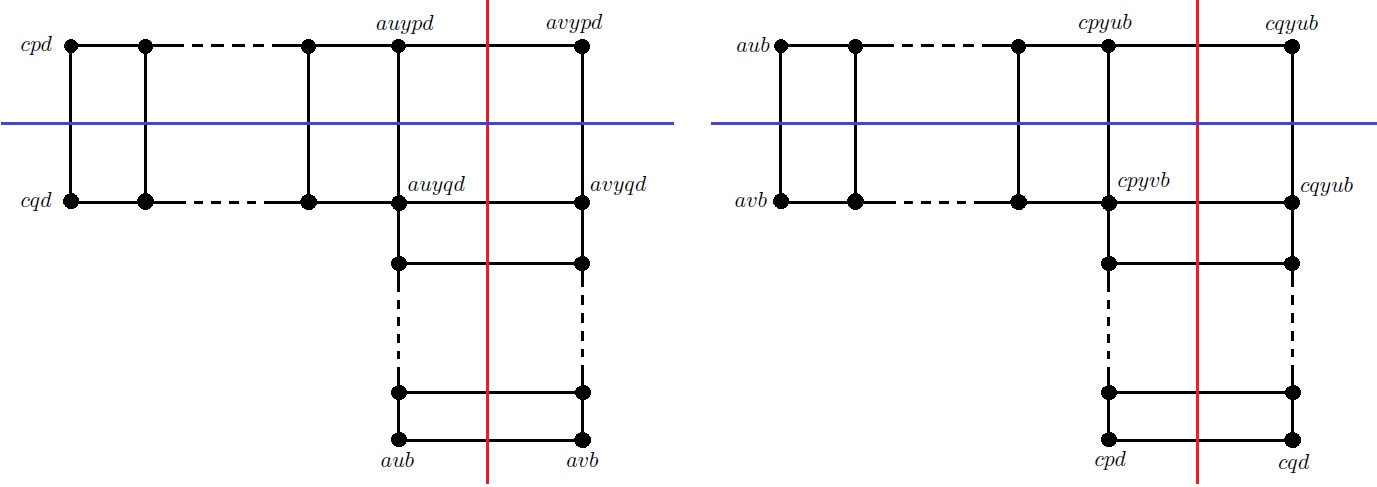}
\end{center}
Conversely, suppose the mentionned hyperplanes are transverse. Two cases happen. If the hyperplanes meet inside a square $(x,p \to q,y,u \to v,z)$ where $[x,p \to q, yuz]=[c,p \to q,d]$ and $[xpy,u \to v,z]=[a,u \to v,b]$, then we deduce from Lemma \ref{hyperplansSquier} that $c=x$, $d=yuz$, $a=xpy$ and $b=z$ modulo $\mathcal{P}$, hence $d=yub$ and $a=cpy$ modulo $\mathcal{P}$. If the hyperplanes meet inside a square $(x,u \to v,y,p \to q,z)$ where $[x,u \to v, ypz]=[a,u \to v,b]$ and $[xuy,p \to q,z]=[c,p \to q,d]$, then we deduce from Lemma \ref{hyperplansSquier} that $a=x$, $b=ypz$, $c=xuy$ and $d=z$ modulo $\mathcal{P}$, hence $b=ypd$ et $c=auy$ modulo $\mathcal{P}$. $\square$

\begin{remark}\label{remarkprec}
Using the relation $\prec$ introduced in Section 4, we have proved more precisely that $[a,u \to v,b] \prec [c,p \to q,d]$ if and only if there exists some word $y \in \Sigma^+$ satisfying $\left\{ \begin{array}{l} c=auy \\ b=ypd \end{array} \right.$ modulo $\mathcal{P}$.
\end{remark}

\noindent
In the sequel, the following notation will be convenient:

\begin{definition}
We denote by $S(\mathcal{P},a)u S(\mathcal{P},b)$ the image in $S(\mathcal{P})$ of the combinatorial map
\begin{center}
$\left\{ \begin{array}{ccc} S(\mathcal{P},a) \times S(\mathcal{P},b) & \to & S(\mathcal{P}) \\ (\alpha, \beta) & \mapsto & \alpha u \beta \end{array} \right.$.
\end{center}
\end{definition}

\noindent
Let $J= [a,u \to v,b]$ be an oriented hyperplane in $S(\mathcal{P},w)$ and let $\tilde{J}$ denote a lift of $J$ in the Farley complex $X(\mathcal{P},w)$. The \textit{neighborhood} $N(\tilde{J})$ of the hyperplane $\tilde{J}$ is defined as the union of the cubes intersecting $\tilde{J}$, and we denote by $\partial \tilde{J}$ the union of the cubes of $N(\tilde{J})$ disjoint from $\tilde{J}$; this subcomplex has two connected components, denoted $\partial_+ \tilde{J}$ and $\partial_- \tilde{J}$ following the natural orientation of $\tilde{J}$ induced by $J$. By extension, we will write $N(J)$, $\partial J$, $\partial_- J$ and $\partial_+ J$ as the images in $S(\mathcal{P},w)$ of $N(\tilde{J})$, $\partial \tilde{J}$, $\partial_- \tilde{J}$ and $\partial_+ \tilde{J}$ respectively.
\begin{center}
\includegraphics[scale=0.6]{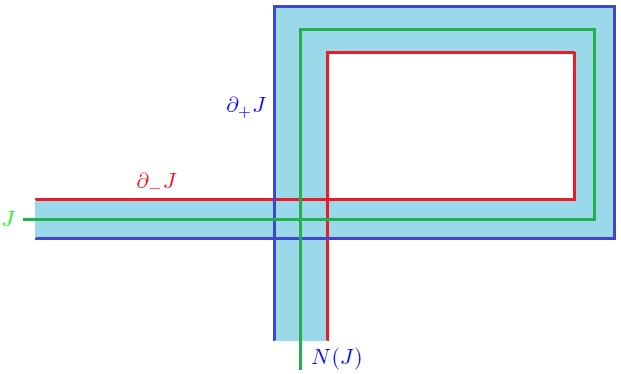}
\end{center}
Our last lemma is a direct consequence of \cite[Lemma 2.5]{arXiv:1505.02053} and \cite[Lemma 1.3.5]{arXiv:1505.02053}. 

\begin{thm}\label{partialJ}
Let $J= [a, u \to v,b]$ be a hyperplane in $S(\mathcal{P})$. Then, $\partial_-J= S(\mathcal{P},a)u S(\mathcal{P},b)$ and $\partial_+J= S(\mathcal{P},a)vS(\mathcal{P},b)$. In particular, if $J$ is a clean hyperplane, then $J$ is naturally isometric to $S(\mathcal{P},a) \times S(\mathcal{P},b)$.
\end{thm}

\subsection{Pathological configurations}

\noindent
A semigroup presentation $\mathcal{P}= \langle \Sigma \mid \mathcal{R} \rangle$ and a base word $w \in \Sigma^+$ being fixed, we determine exactly when pathological configurations of hyperplanes appear in the Squier complex $S(\mathcal{P},w)$. Our main result is:

\begin{thm}\label{Squierspecial}
Let $\mathcal{P}= \langle \Sigma \mid \mathcal{R} \rangle$ be a semigroup presentation and $w_0 \in \Sigma^+$ a base word. Then, the following assertions are equivalent:
\begin{itemize}
	\item[$(i)$] $S(\mathcal{P},w_0)$ is clean,
	\item[$(ii)$] $S(\mathcal{P},w_0)$ has no self-intersecting hyperplanes,
	\item[$(iii)$] there are no words $a,b,p \in \Sigma^+$ such that $w_0=ab$, $a=ap$ and $b=pb$ modulo $\mathcal{P}$ with $[p]_{\mathcal{P}} \neq \{ p\}$.
\end{itemize}
Moreover, $S(\mathcal{P},w_0)$ is special if and only if it satisfies the conditions above and the following one:
\begin{itemize}
	\item[$(iv)$] there are no words $a,u,v,w,b,p,q, \xi \in \Sigma^+$ such that $w_0=auvwb$, $au=au(v \xi)$, $wb=( \xi v)wb$ modulo $\mathcal{P}$ and $uv=p,vw=q \in \mathcal{R}$.
\end{itemize}
\end{thm}

\noindent
A simplified criterion, often sufficient, is the following:

\begin{cor}
If there are no words $a,b,p \in \Sigma^+$ satisfying $w_0=ab$, $a=ap$ and $b=pb$ modulo $\mathcal{P}$, then $S(\mathcal{P},w_0)$ is special. In particular, if $[w_0]_{\mathcal{P}}$ is finite, then $S(\mathcal{P},w_0)$ is special.
\end{cor}

\noindent
\textbf{Proof.} Suppose that $S(\mathcal{P},w_0)$ is not special. Thus, at least one of the points $(iii)$ and $(iv)$ of Theorem \ref{Squierspecial} does not hold. The negation of the point $(iii)$ implies that there exist words $a,b,p \in \Sigma^+$ satisfying $w=ab$, $a=ap$ and $b=pb$ modulo $\mathcal{P}$; and the negation of the point $(iv)$ implies $w_0=(au)(vwb)$, $au=au(v \xi)$ and $vwb=(v \xi) vwb$ modulo $\mathcal{P}$. Therefore, this proves the first assertion of our corollary. 

\medskip \noindent
The second assertion follows from the following observation: if there exist some words $a,b,p \in \Sigma^+$ satisfying $w_0=ab$, $a=ap$ and $b=pb$ modulo $\mathcal{P}$, then $w_0=ap^nb$ modulo $\mathcal{P}$ for every $n \geq 1$, and $[w_0]_{\mathcal{P}}$ is infinite. $\square$

\medskip \noindent
Essentially, Theorem \ref{Squierspecial} will be a consequence of the following lemmas:

\begin{lemma}\label{2sided}
Every hyperplane in $S(\mathcal{P})$ is 2-sided.
\end{lemma}

\noindent
\textbf{Proof.} If there were a 1-sided hyperplane in $S(\mathcal{P})$, there would exist some words $a,x,y,b \in \Sigma^+$ with $[a,x \to y,b]= [a,y \to x,b]$. But such an equality contradicts Lemma \ref{hyperplansSquier}, because $\mathcal{R}$ does not contain the relation $x=x$. $\square$

\begin{lemma}\label{hypautointersectionSquier}
A hyperplane $J$ in $S(\mathcal{P})$ is self-intersecting if and only if $J=[a,p \to q, c]$ for some words $a,p,q,b,c \in \Sigma^+$ satisfying $a=apb$ and $c=bpc$ modulo $\mathcal{P}$.
\end{lemma}

\noindent
\textbf{Proof.} Let $J$ be a hyperplane with $J=[a,p \to q,c]$ for some words $a,b,p,q \in \Sigma^+$ satisfying $a=apb$ and $c=bpc$ modulo $\mathcal{P}$. Then, according to Lemma \ref{hyperplansSquier}, $J$ is dual to the edges $(a,p \to q, bpc)$ and $(apb, p \to q,c)$; since these edges are non-parallel sides of the square $(a,p \to q,b,p \to q,c)$, we deduce that $J$ is self-intersecting. 

\medskip \noindent
Conversely, suppose that $S(\mathcal{P},w)$ contains a self-intersecting hyperplane $J$. Then, there exists a square
\begin{center}
$(a,u \to v,b,p \to q,c)$,
\end{center}
where $[a,u \to v,bpc]=J=[aub,p \to q,c]$. From Lemma \ref{hyperplansSquier}, we deduce that $u=p$ and $v=q$ in $\Sigma^+$, and that $a=aub$ and $c=bpc$ modulo $\mathcal{P}$. Therefore, $J=[a,p \to q,c]$ with the desired equalities modulo $\mathcal{P}$. $\square$

\begin{lemma}\label{hypautotangenceSquier}
A hyperplane $J$ in $S(\mathcal{P})$ is self-osculating if and only if $J=[a,(kh)^nk \to p,b]$ for some $n \geq 1$ and $a,k,h,p,b \in \Sigma^+$ satisfying $a=akh$ and $b=hkb$ modulo $\mathcal{P}$.
\end{lemma}

\noindent
\textbf{Proof.} Let $J$ be a self-osculating hyperplane in $S(\mathcal{P})$. Therefore, $S(\mathcal{P})$ contains the following configuration:
\begin{center}
\includegraphics[scale=0.6]{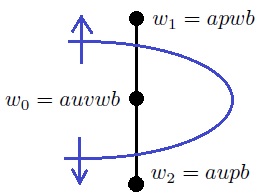}
\end{center}
Indeed, the words $w_1$ and $w_2$ have to be obtained from $w_0$ by applying a same relation of $\mathcal{R}$ on two intersecting subwords. Thus, $uv=vw$ in $\Sigma^+$. Let $n \geq 0$ be the greatest integer such that $u^n$ is a prefix of $v$, ie., $v=u^n k$ for some word $k \in \Sigma^+$. Then, the equality $uv=vw$ becomes $uk=kw$ in $\Sigma^+$. Since $u$ is not a prefix of $k$, by definition of $n$, necessarily $\mathrm{lg}(k)< \mathrm{lg}(u)$, hence $\mathrm{lg}(w)= \mathrm{lg}(u)> \mathrm{lg}(k)$; therefore, $k$ is a suffix of $w$: $w=hk$ for some word $h \in \Sigma^+$. Thus, $u=kh$, $v=(kh)^nk$ and $w=hk$ in $\Sigma^+$, and so $J$ is dual to the edges $(a,(kh)^{n+1}k \to p, hkb)$ and $(akh,(hk)^{n+1}k \to p,b)$. From Lemma \ref{hyperplansSquier}, we deduce that $J=[a,(kh)^nk \to p,b]$ with $a=akh$ and $b=hkb$ modulo $\mathcal{P}$. 

\medskip \noindent
Conversely, let $J$ be a hyperplane such that $J=[a,(kh)^nk \to p,b]$ for some $n \geq 1$ and $a,k,h,p,b \in \Sigma^+$ satisfying $a=akh$ and $b=hkb$ modulo $\mathcal{P}$. In particular, according to Lemma \ref{hyperplansSquier}, $J$ is dual to the edges $(a,(kh)^nk \to p,hkb)$ and $(akh,(kh)^nk \to p,b)$; since these edges have $a(kh)^{n+1}kb$ as a common vertex and do not belong to a same square (because $n \geq 1$), we deduce that $J$ is self-osculating. $\square$

\medskip \noindent
\textbf{Proof of Theorem \ref{Squierspecial}.} The implication $(i) \Rightarrow (ii)$ is clear. 

\medskip \noindent
Now, we prove $(ii) \Rightarrow (iii)$. Suppose that there exist some words $a,b,p \in \Sigma^+$ satisfying $w_0=ab$, $a=ap$ and $b=pb$ modulo $\mathcal{P}$ with $[p]_{\mathcal{P}} \neq \{p \}$. In particular, there exist some words $x,y,r,s \in \Sigma^+$ such that $p=xry$ in $\Sigma^+$ and $r=s \in \mathcal{R}$. Thus, the hyperplane $[ax,r \to s,yb]$ self-intersects inside the square $(ax, r \to s, yx, r \to s,yb)$. A fortiori, the Squier complex $S(\mathcal{P},w)$ contains a self-intersecting hyperplane.

\medskip \noindent
Now, we prove $(iii) \Rightarrow (i)$. The equivalences between $(i)$, $(ii)$ and $(iii)$ will follow.

\medskip \noindent
If $S(\mathcal{P},w_0)$ contains a self-intersecting hyperplane, we deduce from Lemma \ref{hypautointersectionSquier} there exist $a,b,p,c \in \Sigma^+$ such that $w_0=(ap)c$, $ap=(ap)(bp)$ and $c=(bp)c$ modulo $\mathcal{P}$ because $[bp]_{\mathcal{P}} \neq \{bp\}$ (since this class contains $bq$).

\medskip \noindent
If $S(\mathcal{P},w_0)$ contains a self-osculating hyperplane, we deduce from Lemma \ref{hypautointersectionSquier} that there exist $n \geq 1$ and $a,k,h,b \in \Sigma^+$ such that $w_0=akb$, $a=a(kh)^{n+1}$ and $kb=(kh)^{n+1}kb$ modulo $\mathcal{P}$ with $[ (kh)^{n+1}]_{\mathcal{P}} \neq \{ (kh)^{n+1} \}$ (since this class contains $ph$).

\medskip \noindent
Now, we prove that $S(\mathcal{P},w_0)$ is special if and only if conditions $(iii)$ and $(iv)$ are satisfied. 

\medskip \noindent
Suppose that $S(\mathcal{P},w_0)$ is not special. If $S(\mathcal{P},w_0)$ is not clean, we already know that it implies that $(iii)$ is not satisfied. Now, suppose that $S(\mathcal{P},w_0)$ contains two inter-osculating hyperplanes $J_1$ and $J_2$. Two cases may happen: $J_1$ and $J_2$ are respectively dual to edges either of the form $(a,uv \to p,wb)$ and $(au, vw \to q,b)$, or of the form $(au,v \to p, wb)$ and $(a,uvw \to q,b)$. In the first case, we deduce from Lemma \ref{hyptransverseSquier} that there exists a word $\xi$ satisfying either $\left\{ \begin{array}{l} au=auv \xi \\ wb= \xi vw b \end{array} \right.$ modulo $\mathcal{P}$, in which case we have $w_0=(au)(vwb)$, $au=(au)(v\xi)$ and $vwb=(v \xi) vwb$ modulo $\mathcal{P}$ with $uv=p, vw=q \in \mathcal{R}$, so $(iv)$ does not hold; or $\left\{ \begin{array}{l} b= \xi uv wb \\ a=auvw \xi \end{array} \right.$ modulo $\mathcal{P}$, in which case we have $w_0=a(uvwb)$, $a=a(uvw \xi)$ and $uvwb=(uvw \xi)(uvwb)$ modulo $\mathcal{P}$ with $[uvw \xi]_{\mathcal{P}} \neq \{ uvw \xi \}$ (since this class contains $pw \xi$ and $uq \xi$), so $(iii)$ does not hold. In the second case, we deduce from Lemma \ref{hyptransverseSquier} that there exists a word $\xi$ satisfying either $\left\{ \begin{array}{l} a=auv \xi \\ wb= \xi uvwb \end{array} \right.$, in which case we have $w_0=(auv)(wb)$, $auv=(auv)( \xi uv)$ and $wb=(\xi uv)(wb)$ modulo $\mathcal{P}$ with $[\xi uv]_{\mathcal{P}} \neq \{ \xi uv \}$ (since this class contains $ \xi pv$), so $(iii)$ does not hold; or $\left\{ \begin{array}{l} b= \xi vwb \\ au=auvw \xi \end{array} \right.$ modulo $\mathcal{P}$, in which case we have $w_0=(au)(vwb)$, $au=(au)(vw \xi)$ and $vwb=(vw \xi) (vwb)$ modulo $\mathcal{P}$ with $[vw \xi]_{\mathcal{P}} \neq \{ vw \xi \}$ (since this class contains $pw \xi$), so $(iii)$ does not hold. 

\medskip \noindent
Conversely, we already know that, if $(iii)$ does not hold, then $S(\mathcal{P},w_0)$ is not clean, and a fortiori, is not special. Now, suppose that $(iv)$ does not hold, ie., suppose there exist some words $a,u,v,w,b,p,q, \xi \in \Sigma^+$ satisfying $w_0=auvwb$, $au=au(v \xi)$ and $wb= ( \xi v)wb$ modulo $\mathcal{P}$ with $uv=p, vw=q \in \mathcal{R}$. Then, the edges $(a,uv \to p,wb)$ and $(au,vw \to q,b)$ have $auvwb$ as a common endpoint but they do not belong to a same square, and the hyperplanes they define intersect inside the square $(a,uv \to p, \xi , vw \to q,b)$: these two hyperplanes inter-osculate. A fortiori, the Squier complex $S(\mathcal{P},w)$ is not special.  $\square$ 

\medskip \noindent
From the classical theory associated to special cube complexes \cite{MR2377497}, we deduce the following corollaries from Theorem \ref{Squierspecial}:

\begin{cor}\label{RF}
If the following conditions are satisfied:
\begin{itemize}
	\item[$\bullet$] there are no words $a,b,p \in \Sigma^+$ such that $w=ab$, $a=ap$ and $b=pb$ modulo $\mathcal{P}$ with $[p]_{\mathcal{P}} \neq \{ p\}$,
	\item[$\bullet$] there are no words $a,u,v,w,b,p,q, \xi \in \Sigma^+$ such that $w_0=auvwb$, $au=au(v \xi)$, $wb=( \xi v)wb$ modulo $\mathcal{P}$ and $uv=p,vw=q \in \mathcal{R}$,
\end{itemize}
then the diagram group $D(\mathcal{P},w)$ embeds into a right-angled Artin group. In particular, it is linear (and so residually finite).
\end{cor}

\noindent
\textbf{Proof.} The conclusion follows from Theorem \ref{Squierspecial} and Corollary \ref{specialRAAG}. $\square$

\begin{cor}
If $[w]_{\mathcal{P}}$ is finite, then the convex-cocompact subgroups of $D(\mathcal{P},w)$ are separable. In particular, its canonical subgroups are separable.
\end{cor}

\noindent
\textbf{Proof.} The first assertion follows from Theorem \ref{convex-cocompact}. Then, let $H$ be a canonical subgroup, ie., there exist $u_1,\ldots, u_n \in \Sigma^+$ such that $w=u_1 \cdots u_n$ modulo $\mathcal{P}$ and $D(\mathcal{P},u_1) \times \cdots \times D(\mathcal{P},u_n) \simeq H$. Let $\Gamma$ be a $(w,u_1 \cdots u_n)$-diagram. Then
\begin{center}
$(\Delta_1,\ldots, \Delta_n) \mapsto \Gamma \cdot \left( \Delta_1 + \cdots + \Delta_n \right) \cdot \Gamma^{-1}$
\end{center}
defines an isometric embedding $X(\mathcal{P},u_1) \times \cdots \times X(\mathcal{P},u_n) \hookrightarrow X(\mathcal{P},w)$, whose image is a convex subcomplex on which $H$ acts geometrically. Consequently, $H$ is a convex-cocompact subgroup, and so is separable. $\square$

\medskip \noindent
Finally, we are able to deduce the following Tits alternative, since it already holds for (the subgroups of) right-angled Artin groups \cite{MR634562}:

\begin{cor}\label{Titsalternative}
Suppose that the conditions of Corollary \ref{RF} are satisfied. Then any subgroup of $D(\mathcal{P},w)$ is either free abelian or contains a non-abelian free groups. 
\end{cor}

\begin{remark}
Notice that the first point of Corollary \ref{RF} is always satisfied if $S(\mathcal{P},w)$ is finite-dimensional. However, the second point may happen even in this case: if
\begin{center}
$\mathcal{P}=\left\langle \begin{array}{c} a,b,u,v,w, \xi \\ p_1,p_2,p_3, \\ q_1,q_2,q_3 \end{array} \left| \begin{array}{c} au=auv \xi \\ wb= \xi vwb \\ uv=p_1, \ vw=q_1 \end{array}, \ \begin{array}{c} p_1=p_2, p_2=p_3, p_3=p_1 \\ q_1=q_2, q_2=q_3, q_3=q_1 \end{array} \right. \right\rangle$,
\end{center}
then $S(\mathcal{P},auvwb)$ is a 2-dimensional cube complex which is not special, because it contains two inter-osculating hyperplanes. Therefore, Corollary \ref{Titsalternative} does not apply to prove that the associated diagram group 
\begin{center}
$D(\mathcal{P},auvwb) \simeq \langle a,h,t \mid [t,a^{h^n}]=1, \ n \geq 1 \rangle$
\end{center}
satisfies the Tits alternative. In fact, using almost verbatim \cite[Section 4]{MR1703363} combined with the relation $\prec$ we introduce in Section 4, it can be proved that the conclusion of Corollary \ref{Titsalternative} holds whenever $S(\mathcal{P},w)$ is finite-dimensional.
\end{remark}

\subsection{Explicit embedding}

\noindent
Corollary \ref{RF} proves that some diagram groups embed into a right-angled Artin group; more precisely, it proves that the morphism of Theorem \ref{localisometry} is injective. In this section, we show how to describe this morphism explicitly. That is to say, given a precise example of a semigroup presentation $\mathcal{P}= \langle \Sigma \mid \mathcal{R} \rangle$ and a base word $w \in \Sigma^+$ (such that the Squier complex $S(\mathcal{P},w)$ is special), we want to be able to draw a graph $\Gamma$, such that there exists an embedding $D(\mathcal{P},w) \hookrightarrow A(\Gamma)$, and to write down the images in $A(\Gamma)$ of some generating set of $D(\mathcal{P},w)$.

\begin{definition}\label{Phi}
Let $A(\mathcal{P},w)$ be the right-angled Artin group associated to the graph $\Gamma(\mathcal{P},w)$ whose vertices are the (unoriented) hyperplanes of the Squier complex $S(\mathcal{P},w)$ and whose edges link two distinct intersecting hyperplanes.

\medskip \noindent
When $S(\mathcal{P},w)$ has no self-intersecting hyperplanes, $\Gamma(\mathcal{P},w)$ is the transversality graph of the Squier complex $S(\mathcal{P},w)$.

\medskip \noindent
Fix an orientation of the edges of $S(\mathcal{P},w)$. Then, the map sending each positive edge $(a,u \to v,b)$ of $S(\mathcal{P},w)$ to $[a,u \to v,b] \in A(\mathcal{P},w)$, and each negative edge $(a,u \to v,b)$ of $S(\mathcal{P},w)$ to $[a,u \to v,b]^{-1} \in A(\mathcal{P},w)$, induces a morphism from $\pi_1( S(\mathcal{P},w),w)$ to $A(\mathcal{P},w)$, and so a morphism
\begin{center}
$\Phi= \Phi(\mathcal{P},w) : D(\mathcal{P},w) \to A(\mathcal{P},w)$.
\end{center}
Notice that the isomorphism between the fundamental group $\pi_1( S(\mathcal{P},w),w)$ and the diagram group $D(\mathcal{P},w)$ is made explicit by the discussion preceding Lemma \ref{iso}, so that $\Phi$ can be describe explicitly.
\end{definition}

\noindent
In fact, when $S(\mathcal{P},w)$ is special, the morphism $\Phi$ is exactly the one used to prove Theorem \ref{localisometry}, so

\begin{prop}\label{Phispecial}
If $S(\mathcal{P},w)$ is special, then $\Phi : D(\mathcal{P},w) \to A(\mathcal{P},w)$ is injective.
\end{prop}

\begin{ex}\label{ZbulletZspecial}
Let us consider the following semigroup presentation:
\begin{center}
$\mathcal{P} = \left\langle \begin{array}{l} a_1,a_2,a_3 \\ b_1,b_2,b_3 \end{array}, p \left| \begin{array}{l} a_1=a_2, a_2=a_3, a_3=a_1 \\ b_1=b_2, b_2=b_3, b_3=b_1 \end{array}, \begin{array}{l} a_1=a_1p \\ b_1=pb_1 \end{array} \right. \right\rangle$. 
\end{center}
The diagram group $D(\mathcal{P},a_1b_1)$ is denoted by $\mathbb{Z} \bullet \mathbb{Z}$ in \cite[Section 8]{MR1396957}; it is a group which is finitely-generated but not finitely-presented, with 
\begin{center}
$\langle a,b,z \mid [a,b^{z^n}]=1, \ n \geq 0 \rangle$ 
\end{center}
as a presentation (see \cite[Lemma 8.5]{MR1396957} or Example \ref{ZbulletZgraph}).

\medskip \noindent
According to Theorem \ref{Squierspecial}, the Squier complex $S(\mathcal{P},a_1b_1)$ is special, so that $\mathbb{Z} \bullet \mathbb{Z}$ is embeddable into a right-angled Artin. Using the discussion above, now we want to describe explicitly such an embedding.

\medskip \noindent
Using Lemma \ref{hyperplansSquier}, we find that $S(\mathcal{P},a_1b_1)$ has eight hyperplanes: $A_i=[1,a_i \to a_{i+1},b_1]$, $B_i=[a_1,b_i \to b_{i+1},1]$, $C= [ 1,a_1 \to a_1p,b_1]$ and $D=[a_1, b_1 \to pb_1,1]$. Using Lemma \ref{hyptransverseSquier}, we find that $\Gamma(\mathcal{P},a_1b_1)$ is a complete bipartite graph $K_{4,4}$, where each vertex of $\{ A_1,A_2,A_3,C \}$ is linked by an edge to each vertex of $\{ B_1, B_2, B_3, D \}$. In particular, $A(\mathcal{P},a_1b_1) \simeq \mathbb{F}_4 \times \mathbb{F}_4$. 

\medskip \noindent
Then, using \cite[Theorem 9.8]{MR1396957} or Example \ref{ZbulletZgraph}, we find that $\mathbb{Z} \bullet \mathbb{Z}$ is generated by the three following diagrams:
\begin{center}
\includegraphics[scale=0.6]{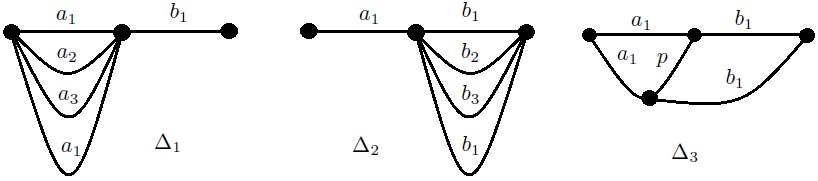}
\end{center}
For instance, according to Lemma \ref{iso}, the first diagram corresponds to loop of edges:
$(1,a_1 \to a_2, b_1), \ (1,a_2 \to a_3, b_1), \ (1,a_3 \to a_1, b_1)$,
so that $\Phi(\Delta_1)=A_1A_2A_3^{-1}$. In the same way, we find that $\Phi(\Delta_2)= B_1B_2B_3^{-1}$ and $\Phi(\Delta_3)=CD^{-1}$.

\medskip \noindent
We conclude that the subgroup $\langle A_1A_2A_3^{-1}, B_1B_2B_3^{-1}, CD^{-1} \rangle$ of 
\begin{center}
$\mathbb{F}_4 \times \mathbb{F}_4= \langle A_1,A_2,A_3,C \mid \ \rangle \ast \langle B_1,B_2,B_3, D \mid \ \rangle$
\end{center}
is isomorphic to $\mathbb{Z} \bullet \mathbb{Z}$.

\medskip \noindent
However, it is worth noticing that, by collapsing $\{ A_1, A_2, A_3 \}$ and $\{ B_1, B_2, B_3 \}$ respectively to the points $A$ and $B$, the embedding $\Phi : \mathbb{Z} \bullet \mathbb{Z} \hookrightarrow \mathbb{F}_4 \times \mathbb{F}_4$ may be simplified into:
\begin{center}
\includegraphics[scale=0.6]{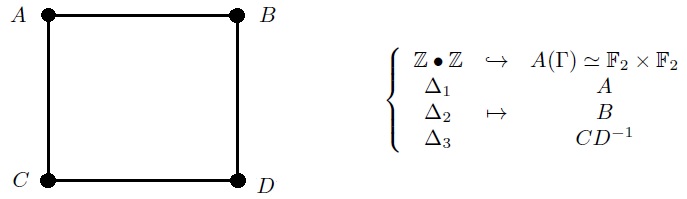}
\end{center}
Therefore, if
\begin{center}
$\mathbb{F}_2 \times \mathbb{F}_2= \langle a,b \mid \ \rangle \ast \langle x,y | \ \rangle = \langle x,y,a,b \mid [x,a]=[x,b]=[y,a]=[y,b]=1 \rangle$,
\end{center}
then the subgroup $\langle a,x,by \rangle$ is isomorphic to $\mathbb{Z} \bullet \mathbb{Z}$. 
\end{ex}

\section{Quasi-isometric embeddability into a product of trees}

\subsection{Rank of a hyperplane}

\noindent
In this Section, we fix a semigroup presentation $\mathcal{P}= \langle \Sigma \mid \mathcal{R} \rangle$ and a base word $w \in \Sigma^+$.

\begin{definition}
Let $J_1$ and $J_2$ be two hyperplanes in the Squier complex $S(\mathcal{P},w)$. If they meet inside a square $(a,u \to v, b,p \to q,c)$ so that $[a,u \to v,bpc]=J_1$ and $[aub,p \to q,c]=J_2$, we write $J_1 \prec J_2$.
\end{definition}

\noindent
Remark \ref{remarkprec} gives a necessary and sufficient condition for two hyperplanes $J_1$ and $J_2$ to satisfy $J_1 \prec J_2$.

\begin{lemma}\label{lemprec}
Let $J_1, \ldots, J_n$ be $n$ hyperplanes satisfying $J_i \prec J_{i+1}$ for all $1 \leq i \leq n-1$. Then there exists an $n$-cube
\begin{center}
$(a_1,u_1 \to v_1, a_2, u_2 \to v_2, \ldots, a_n, u_n \to v_n, a_{n+1})$
\end{center}
such that $[a_1u_1 \cdots a_{k-1}u_{k-1}a_k, u_k \to v_k, a_{k+1}u_{k+1} \cdots a_nu_na_{n+1}]=J_k$. 
\end{lemma}

\noindent
\textbf{Proof.} We prove the lemma by induction on $n$. By definition of $\prec$, the result is true for $n=2$. Now, let $J_1, \ldots, J_{n+1}$ be $n+1$ hyperplanes satisfying $J_i \prec J_{i+1}$ for all $1 \leq i \leq n$. By our induction hypothesis, there exists an $n$-cube
\begin{center}
$(a_1,u_1 \to v_1, a_2, u_2 \to v_2, \ldots, a_n, u_n \to v_n, a_{n+1})$
\end{center}
such that $[a_1u_1 \cdots a_{k-1}u_{k-1}a_k, u_k \to v_k, a_{k+1}u_{k+1} \cdots a_nu_na_{n+1}]=J_k$ for all $1 \leq k \leq n$. Then, because $J_n \prec J_{n+1}$, there exist a square
\begin{center}
$(a,u_n \to v_n, b, p \to q,c)$
\end{center}
such that $[a, u_n \to v_n,bpc]=J_n$ and $[au_nb,p \to q,c]=J_{n+1}$. Because we also have $J_n=[a_1u_1 \cdots a_{n-1}u_{n-1}a_n,u_n \to v_n, a_{n+1}]$, we deduce from Lemma \ref{hyperplansSquier} that $a=a_1u_1 \cdots a_{n-1}u_{n-1}a_n$ and $bpc=a_{n+1}$ modulo $\mathcal{P}$. This proves the existence of the $(n+1)$-cube
\begin{center}
$(a_1,u_1 \to v_1, a_2, u_2 \to v_2, \ldots, a_n, u_n \to v_n, b,p \to q,c)$
\end{center}
Noticing that 
\begin{center}
$\begin{array}{l} [a_1u_1 \cdots a_{k-1}u_{k-1}a_k, \ u_k \to v_k, \ a_{k+1}u_{k+1} \cdots a_nu_nbpc] \\ = [a_1u_1 \cdots a_{k-1}u_{k-1}a_k, \ u_k \to v_k, \ a_{k+1}u_{k+1} \cdots a_nu_na_{n+1}]= J_k \end{array}$ 
\end{center}
and $[a_1u_1 \cdots a_nu_nb,p \to q,c]= [au_nb,p \to q,c]=J_{n+1}$ proves our lemma. $\square$

\begin{cor}\label{propprec}
The relation $\prec$ satisfies the following properties:
\begin{itemize}
	\item[$\bullet$] If $J_1 \prec J_2$ and $J_2 \prec J_3$, then $J_1 \prec J_3$.
	\item[$\bullet$] $J_1$ and $J_2$ are comparable with respect to $\prec$ if and only if they intersect. 
	\item[$\bullet$] $\mathrm{max} \{ n \geq 0 \mid \text{there exist} \ J_1, \ldots, J_n \ \text{such that} \ J_1 \prec \cdots \prec J_n \}= \dim S(\mathcal{P},w)$.
\end{itemize}
\end{cor}

\noindent
\textbf{Proof.} Let $J_1,J_2,J_3$ be three hyperplanes such that $J_1 \prec J_2$ and $J_2 \prec J_3$. According to Lemma \ref{lemprec}, there exists a cube
\begin{center}
$(a,u \to v, b, p \to q, c, x \to y,d)$
\end{center}
such that $[a,u \to v,bpcxd]=J_1$, $[aub,p \to q, cxd]=J_2$ and $[aubpc,x \to y,d]=J_3$. By considering the square
\begin{center}
$(a,u \to v, bpc, x \to y,d)$,
\end{center}
we prove the first point.

\medskip \noindent
The second point is clear by definition of $\prec$. 

\medskip \noindent
Let us prove the third point. If there exist $n$ hyperplanes $J_1, \ldots, J_n$ satisfying $J_1 \prec \cdots \prec J_n$, Lemma \ref{lemprec} yields an $n$-cube, hence $\dim S(\mathcal{P},w) \geq n$. If $n\leq \dim S(\mathcal{P},w)$, there exists a $n$-cube
\begin{center}
$(a_1,u_1 \to v_1,a_2, u_2 \to v_2,\ldots, a_n,u_n \to v_n,a_{n+1})$
\end{center}
in $S(\mathcal{P},w)$. Let $J_k= [ a_1u_1 \cdots a_{k-1}u_{k-1}a_k, u_k \to v_k, a_{k+1}u_{k+1} \cdots a_nu_na_{n+1}]$. Then, it is clear that $J_1 \prec \cdots \prec J_n$. Therefore, we have proved the equality
\begin{center}
$\mathrm{max} \{ n \geq 0 \mid \text{there exist} \ J_1, \ldots, J_n \ \text{such that} \ J_1 \prec \cdots \prec J_n \}= \dim S(\mathcal{P},w)$.
\end{center}
The proof is complete. $\square$

\medskip \noindent
It is worth noticing that the relation $\prec$ gives restrictions on the geometry of Squier complexes, namely on their \textit{transversality graphs}.

\begin{definition}
The \textit{transversality graph} of a cube complex $X$ is defined as the graph whose vertices are the hyperplanes of $X$ and whose edges link two transverse hyperplanes.
\end{definition}

\begin{cor}\label{cycle}
The transverality graph of $S(\mathcal{P},w)$ has no induced cycles of odd length greater than five.
\end{cor}

\noindent
Recall that $\Gamma \subset \Lambda$ is an \textit{induced subgraph} of $\Lambda$ if any vertices $x,y \in \Gamma$ are linked by an edge in $\Lambda$ if and only if they are linked by an edge in $\Gamma$.

\medskip \noindent
\textbf{Proof.} Let $n \geq 5$ be an odd integer. Suppose by contradiction that there exist $n$ hyperplanes $J_1, \ldots, J_n$ such that $J_i$ and $J_k$ are transverse if and only if $k=i \pm 1$ (modulo $n$). Suppose that $J_1 \prec J_2$; the case $J_2 \prec J_1$ will be completely symmetric.

\medskip \noindent
Since $J_2$ and $J_3$ are transverse, either $J_2 \prec J_3$ or $J_3 \prec J_2$. But we already know that $J_1 \prec J_2$ so that $J_2 \prec J_3$ would imply $J_1 \prec J_3$ and a fortiori that $J_1$ and $J_3$ are transverse. Therefore, $J_2 \prec J_3$. Similarly, we deduce that $J_2 \prec J_4$, $J_5 \prec J_4$, and so on. Thus, $J_{2k+1} \prec J_{2k}$ for all $0 \leq k \leq \frac{n-1}{2}$. In particular, $J_n \prec J_{n-1}$ since $n$ is odd. 

\medskip \noindent
Then, because $J_n$ and $J_1$ are transverse, either $J_1 \prec J_n$ or $J_n \prec J_1$. In the first case, we deduce from $J_n \prec J_{n-1}$ that $J_1$ and $J_{n-1}$ are transverse, a contradiction. In the second case, we deduce from $J_1 \prec J_2$ that $J_n$ and $J_2$ are transverse, a contradiction. $\square$

\begin{remark}
Corollary \ref{cycle} does not hold for induced cycles of even length. For every even integer $n \geq 2$, it is possible to find a Squier complex containing an induced cycle of length $n$.
\end{remark}

\noindent
At least in the finite-dimensional case, the relation $\prec$ allows to distinguish families of hyperplanes according to their \textit{ranks}:

\begin{definition}
Let $J$ be a hyperplane in the Squier complex $S(\mathcal{P},w)$. We define its \textit{rank} as
\begin{center}
$\mathrm{rank}(J) = \sup \{ n \geq 0 \mid \text{there exist} \ J_1, \ldots, J_n \ \text{with} \ J_1 \prec \cdots \prec J_n \prec J \}$.
\end{center}
By extension, we define the \textit{rank} of a hyperplane in the Farley complex $X(\mathcal{P},w)$ by the rank of its image by the covering map $X(\mathcal{P},w) \to S(\mathcal{P},w)$.
\end{definition}

\noindent
Mainly, we will be interested in the case where the rank of any hyperplane is finite, ie., when there does not exist an infinite sequence $\cdots \prec J_2 \prec J_1$. For instance, it happens when $S(\mathcal{P},w)$ is finite-dimensional, or more generally, when $S(\mathcal{P},w)$ does not contain an infinite family of pairwise intersecting hyperplanes. Notice also that the rank of a self-intersecting hyperplane is always infinite, so that our cube complexes will be always clean when we will want the hyperplanes' ranks to be well-defined.

\begin{lemma}
Let $J_1$ and $J_2$ be two hyperplanes in the Farley complex $X(\mathcal{P},w)$ of finite rank. If $\mathrm{rank}(J_1)= \mathrm{rank}(J_2)$ then $J_1$ and $J_2$ are disjoint.
\end{lemma}

\noindent
\textbf{Proof.} Suppose that $J_1$ and $J_2$ are transverse. Then they are comparable with respect to $\prec$, say $J_1 \prec J_2$. If $r= \mathrm{rank}(J_1)$, let $H_1,\ldots, H_r$ be $r$ hyperplanes such that $H_1 \prec \cdots \prec H_r \prec J_1$. Then,
\begin{center}
$H_1 \prec \cdots \prec H_r \prec J_1 \prec J_2$,
\end{center}
hence $\mathrm{rank}(J_2)> \mathrm{rank}(J_1)$. In particular, $\mathrm{rank}(J_2) \neq \mathrm{rank}(J_1)$. $\square$

\medskip \noindent
Let $\mathfrak{J}$ denote the set of hyperplanes of the Farley complex $X(\mathcal{P},w)$, and, for every $k \geq 0$, let $\mathfrak{J}_k$ denote the subset of hyperplanes of rank $k$. As a consequence of the previous lemma, we deduce that $\mathfrak{J}_k$ induces an \textit{arboreal structure} on $X(\mathcal{P},w)$, ie., the graph whose vertices are the connected components of $X(\mathcal{P},w) \backslash \mathfrak{J}_k$ and whose edges link two adjacent components is a tree. Equivalently, the cube complex constructed by cubulation from the pocset defines by the halfspaces delimited by the hyperplanes of $\mathfrak{J}_k$ is a tree $\Lambda_k$. 

\medskip \noindent
There is a natural map $X(\mathcal{P},w)^{(0)} \to \Lambda_k^{(0)}$, sending every vertex of $X(\mathcal{P},w)$ to the principal ultrafilter it defines, which extends to a combinatorial map $\Phi_k : X(\mathcal{P},w) \to \Lambda_k$. 

\newpage

\noindent
Let 
\begin{center}
$\Phi = \Phi_0 \times \Phi_1 \times \cdots  : X(\mathcal{P},w) \to \Lambda_0 \times \Lambda_1 \times \cdots$
\end{center}
be the product of these maps.

\subsection{Property B}

\noindent
We define below Property B, as introduced in \cite{MR2271228}. 

\begin{definition}
Let $D(\mathcal{P},w)$ be a diagram group with a finite generating set $S$; $| \cdot|$ will denote the word length function associated to $S$ and $\#( \cdot)$ the number of cells of a semigroup diagram. We say that $D(\mathcal{P},w)$ satisfies \textit{Property B} if there exist $C_1,C_2>0$ such that
\begin{center}
$C_1 \cdot \#( \Delta) \leq | \Delta | \leq C_2 \cdot \# ( \Delta )$
\end{center}
for all spherical diagram $\Delta \in D(\mathcal{P},w)$. 
\end{definition}

\noindent
Property B is used in \cite{MR2271228} to compute the \textit{equivariant uniform Hilbert space compression} of some diagram groups.

\begin{definition}
Let $G$ be a finitely-generated group with $| \cdot |$ a word length function associated to some finite generating set. The \textit{uniform Hilbert space compression} of $G$ is defined as the supremum of the $\alpha$'s such that there exist a Hilbert space $\mathcal{H}$, an embedding $f : G \to \mathcal{H}$ and some constants $C_1,C_2>0$ satisfying
\begin{center}
$C_1 \cdot d(x,y)^{\alpha} \leq \| f(x)-f(y) \| \leq C_2 \cdot d(x,y)$
\end{center}
for all $x,y \in G$.

\medskip \noindent
Similarly, the \textit{equivariant uniform Hilbert space compression} of $G$ is defined by requiring $f$ to be furthermore $G$-invariant.
\end{definition}

\begin{thm}\label{AGS}
\emph{\cite[Theorem 1.13]{MR2271228}} The equivariant uniform Hilbert space compression of a finitely generated diagram group with Property B is a least $1/2$.
\end{thm}

\noindent
In \cite{MR2271228}, it is proved that the Thompson's group $F$ and the lamplighter group $\mathbb{Z} \wr \mathbb{Z}$ satisfy Property B. In fact, the following problem is mentionned \cite[Question 1.6]{MR2271228}:

\begin{question}
Do all finitely generated diagram groups satisfy Property B?
\end{question}

\noindent
Below, we give an equivalent characterization of Property B, so that we will be able to give an alternative proof of Theorem \ref{AGS} and new examples of diagram groups satisfying Property B.

\begin{lemma}\label{propB}
Let $\mathcal{P}= \langle \Sigma \mid \mathcal{R} \rangle$ be a semigroup presentation and $w \in \Sigma^+$ a base word. Suppose that the diagram group $D(\mathcal{P},w)$ is finitely-generated. Then $D(\mathcal{P},w)$ satisfies Property B if and only if the canonical map $D(\mathcal{P},w) \to X(\mathcal{P},w)$, sending a spherical diagram to the vertex of $X(\mathcal{P},w)$ it defines, is a quasi-isometric embedding with respect to the combinatorial metric. 
\end{lemma}

\noindent
\textbf{Proof.} According to \cite[Corollary 2.4]{arXiv:1505.02053},
\begin{center}
$(A,B) \mapsto \#(A^{-1} \cdot B)$,
\end{center}
where $A,B \in X(\mathcal{P},w)$, coincides with the combinatorial distance on $X(\mathcal{P},w)$, so that the conclusion follows. $\square$ 

\begin{cor}\label{cocompactpropB}
Let $\mathcal{P}= \langle \Sigma \mid \mathcal{R} \rangle$ be a semigroup presentation and $w \in \Sigma^+$ a base word. If $[w]_{\mathcal{P}}$ is finite, then $D(\mathcal{P},w)$ satisfies Property B.
\end{cor}

\noindent
\textbf{Proof.} Since $[w]_{\mathcal{P}}$ is finite, the action of $D(\mathcal{P},w)$ on $X(\mathcal{P},w)$ is properly discontinuous and cocompact (Proposition \ref{Farleyaction}). Therefore, according to Milnor-\v Svarc lemma, the map
\begin{center}
$\left\{ \begin{array}{ccc} D(\mathcal{P},w) & \to & X(\mathcal{P},w) \\ \Delta & \mapsto & \Delta \cdot \epsilon(w) \end{array} \right.$
\end{center}
is a quasi-isometry. But it coincides with the canonical map $D(\mathcal{P},w) \to X(\mathcal{P},w)$, so that $D(\mathcal{P},w)$ satisfies Property B according to Lemma \ref{propB}. $\square$

\medskip \noindent
Below is the sketch of an alternative proof of Theorem \ref{AGS}.

\medskip \noindent
\textbf{Proof of Theorem \ref{AGS}.} Let $\mathfrak{J}$ denote the set of hyperplanes in the Farley complex $X(\mathcal{P},w)$. For every vertex $\Delta \in X(\mathcal{P},w)$, we define the map 
\begin{center}
$w_{\Delta} : \left\{ \begin{array}{ccc} \mathfrak{J} & \to & \{ 0,1 \} \\ J & \mapsto & \left\{ \begin{array}{cl} 1 & \text{if} \ J \ \text{separates} \ \Delta \ \text{and} \ \epsilon(w) \\ 0 & \text{otherwise} \end{array} \right. \end{array} \right.$.
\end{center}
Now, following \cite{MR1459140}, the map
\begin{center}
$f : \left\{ \begin{array}{ccc} X(\mathcal{P},w)^{(0)} & \to & \ell^2(\mathfrak{J}) \\ \Delta & \mapsto & \sum\limits_{J \in \mathfrak{J}} w_{\Delta}(J) \cdot \delta_J \end{array} \right.$,
\end{center}
where $\delta_J : H \mapsto \left\{ \begin{array}{cl} 1 & \text{if} \ J=H \\ 0 & \text{otherwise} \end{array} \right.$, is $D(\mathcal{P},w)$-invariant and satisfies
\begin{center}
$\| f(x)-f(y) \|_{\ell^2(\mathfrak{J})} = \sqrt{d_c(x,y)}$
\end{center}
for all $x,y \in X(\mathcal{P},w)^{(0)}$. Now, since $D(\mathcal{P},w)$ satisfies Property B, from Lemma \ref{propB} we deduce that this group quasi-isometrically embeds into $X(\mathcal{P},w)^{(0)}$ with respect to the combinatorial distance $d_c$, so that $f$ induces an equivariant uniform embedding $D(\mathcal{P},w) \hookrightarrow \ell^2(\mathfrak{J})$ whose associated compression is $1/2$. $\square$

\medskip \noindent
We conclude this section with a last example of a finitely-generated diagram group satisfying Property B.

\begin{lemma}\label{ZbulletZpropertyB}
$\mathbb{Z} \bullet \mathbb{Z}$ satisfies Property B.
\end{lemma}

\noindent
Here, $\mathbb{Z} \bullet \mathbb{Z}$ is canonically interpreted as the diagram group given in Example \ref{ZbulletZspecial}. In particular, the class of the base word modulo the semigroup presentation is infinite, so that Corollary \ref{cocompactpropB} does not apply.

\medskip \noindent
\textbf{Sketch of proof.} We use the normal form associated to the decomposition of $\mathbb{Z} \bullet \mathbb{Z}$ as an HNN extension given in Example \ref{ZbulletZgraph} to prove that
\begin{center}
$\frac{1}{2} | \cdot | \leq \# (\cdot) \leq 3 | \cdot |$,
\end{center}
where $| \cdot |$ is the word length function associated to the generating set $\{a,h,t \}$ associated to the presentation
\begin{center}
$\mathbb{Z} \bullet \mathbb{Z} = \langle a,h,t \mid [a,h^{t^n}]=1, \ n \geq 0 \rangle$
\end{center}
we mentionned in Example \ref{ZbulletZspecial}. $\square$

\subsection{Embeddings into a product of trees}

\begin{thm}\label{Phitrees}
Suppose there does not exist any infinite descending chain of hyperplanes of $S(\mathcal{P},w)$ with respect to $\prec$. Then, the combinatorial map 
\begin{center}
$\Phi : X(\mathcal{P},w) \to \Lambda_0 \times \Lambda_1 \times \cdots$ 
\end{center}
is an isometric embedding with respect to the combinatorial metrics.
\end{thm}

\noindent
\textbf{Proof.} The set of hyperplanes of $\Lambda= \Lambda_0 \times \cdots$ may be written as the disjoint union $\mathfrak{H}_0 \sqcup \mathfrak{H}_1 \sqcup \cdots$, where $\mathfrak{H}_k$ corresponds to the set of hyperplanes transverse to the factor $\Lambda_k$. Because the combinatorial distance corresponds to the number of hyperplanes separating two given vertices, we deduce that
\begin{center}
$\displaystyle d_{\Lambda}(x,y) = \sum\limits_{k\geq 0} \# \{ J \in \mathfrak{H}_k \ \text{separating} \ x \ \text{and} \ y \}$
\end{center}
for all $x,y \in \Lambda$. On the other hand, according to Lemma \ref{distancequotient}, for all $x,y \in X(\mathcal{P},w)$, 
\begin{center}
$d_{\Lambda_k}(\Phi_k(x),\Phi_k(y))= \# \{ J \in \mathfrak{J}_k \ \text{separating} \ x \ \text{and} \ y \}$.
\end{center}
However, there is a natural bijection $\mathfrak{J}_k \to \mathfrak{H}_k$, so that
\begin{center}
$d_{\Lambda}(\Phi(x),\Phi(y))=d(x,y)$
\end{center}
for all $x,y \in X(\mathcal{P},w)$. $\square$

\begin{remark}\label{Zinfty}
Theorem \ref{Phitrees} always applies when the Squier complex $S(\mathcal{P},w)$ is finite-dimensional, because of Corollary \ref{propprec}. However, there are also interesting infinite-dimensional examples when it applies. For instance, let us consider the semigroup presentation
\begin{center}
$\mathcal{P} = \langle x,a,b,c \mid x=xa, a=b, b=c, c=a \rangle$.
\end{center}
Then the Squier complex $S(\mathcal{P},x)$ is infinite-dimensional (notice that $x=xa^n$ modulo $\mathcal{P}$, with $D(\mathcal{P},a) \neq \{1 \}$, and apply Proposition \ref{dimSquier} below), but it does not contain any infinite descending chain of hyperplanes with respect to $\prec$ (use Remark \ref{remarkprec}). The diagram group $D(\mathcal{P},x)$ is the free abelian group $\mathbb{Z}^{\infty}$ of infinite (countable) rank.
\end{remark}

\noindent
From the previous theorem, we now deduce the following result:

\begin{thm}\label{QItrees}
Let $\mathcal{P}= \langle \Sigma \mid \mathcal{R} \rangle$ be a semigroup presentation and $w \in \Sigma^+$ a base word. Suppose that $S(\mathcal{P},w)$ is finite-dimensional and $D(\mathcal{P},w)$ finitely-generated. If $D(\mathcal{P},w)$ satisfies Property B, then it quasi-isometrically embeds into a product of $\dim S(\mathcal{P},w)$ trees.
\end{thm}

\noindent
\textbf{Proof.} Because $D(\mathcal{P},w)$ satisfies Property B, it quasi-isometrically embeds into $X(\mathcal{P},w)$ with respect to the combinatorial distance (Lemma \ref{propB}). Therefore, the conclusion follows from Theorem \ref{Phitrees}. $\square$

\medskip \noindent
As direct consequences of Theorem \ref{QItrees}, we have:

\begin{cor}
Let $\mathcal{P}= \langle \Sigma \mid \mathcal{R} \rangle$ be a semigroup presentation and $w \in \Sigma^+$ a base word. Suppose that $S(\mathcal{P},w)$ is finite-dimensional and $D(\mathcal{P},w)$ finitely-generated. If $D(\mathcal{P},w)$ satisfies Property B, then its asymptotic dimension is bounded above by $\dim S(\mathcal{P},w)$.
\end{cor}

\begin{cor}
Let $\mathcal{P}= \langle \Sigma \mid \mathcal{R} \rangle$ be a semigroup presentation and $w \in \Sigma^+$ a base word. Suppose that $S(\mathcal{P},w)$ is finite-dimensional and $D(\mathcal{P},w)$ finitely-generated. If $D(\mathcal{P},w)$ satisfies Property B, then its uniform Hilbert space compression is $1$.
\end{cor}

\noindent
We conclude this section with the following remark. Theorem \ref{QItrees} applies when the Squier complex $S(\mathcal{P},w)$ is finite-dimensional, so a natural question is: given $\mathcal{P}$ and $w$, how to determine whether or not $S(\mathcal{P},w)$ is finite-dimensional? And if it is the case, what is its dimension? A simple criterion is given by the proposition below.

\begin{prop}\label{dimSquier}
Let $n \geq 1$. Then $S(\mathcal{P},w)$ has dimension at least $n$ if and only if there exist some words $u_1,\ldots, u_n \in  \Sigma^+$ such that $w=u_1 \cdots u_n$ modulo $\mathcal{P}$ with $[u_i]_{\mathcal{P}} \neq \{ u_i \}$ for all $1 \leq i \leq n$.
\end{prop}

\noindent
\textbf{Proof.} If $\dim S(\mathcal{P},w) \geq n$, then $S(\mathcal{P},w)$ contains an $n$-cube
\begin{center}
$(a_1,p_1 \to q_1, \ldots, a_n,p_n \to q_n,a_{n+1})$.
\end{center}
Then, $w=(a_1p_1) \cdots (a_{n-1}p_{n-1})(a_np_na_{n+1})$ modulo $\mathcal{P}$ with $[a_ip_i]_{\mathcal{P}} \neq \{ a_ip_i \}$ since $a_iq_i \in [a_ip_i]_{\mathcal{P}}$ and $[a_np_na_{n+1}]_{\mathcal{P}} \neq \{ a_np_na_{n+1} \}$ since $a_nq_na_{n+1} \in [a_np_na_{n+1}]_{\mathcal{P}}$. 

\medskip \noindent
Conversely, suppose there exist some words $u_1,\ldots, u_n \in  \Sigma^+$ such that $w=u_1 \cdots u_n$ modulo $\mathcal{P}$ with $[u_i]_{\mathcal{P}} \neq \{ u_i \}$ for all $1 \leq i \leq n$. For each $1 \leq i \leq n$, there exist $x_i,y_i,p_i,q_i \in \Sigma^+$ such that $u_i=x_ip_iy_i$ in $\Sigma^+$ and $p_i=q_i \in \mathcal{R}$. Then
\begin{center}
$(x_1,p_1 \to q_1,y_1x_2,p_2 \to q_2, \ldots, y_{n-1}x_n,p_n \to q_n,y_n)$
\end{center}
defines an $n$-cube in $S(\mathcal{P},w)$, hence $\dim S(\mathcal{P},w) \geq n$. $\square$

\begin{cor}\label{ZbulletZtrees}
$\mathbb{Z} \bullet \mathbb{Z}$ quasi-isometrically embeds into a product of two trees.
\end{cor}

\noindent
\textbf{Proof.} We interprete $\mathbb{Z} \bullet \mathbb{Z}$ as the diagram group $D(\mathcal{P},a_1b_1)$ given by Example \ref{ZbulletZspecial}. According to Lemma \ref{ZbulletZpropertyB} and Theorem \ref{QItrees}, $\mathbb{Z} \bullet \mathbb{Z}$ quasi-isometrically embeds into a product of $\dim S(\mathcal{P},a_1b_1)$ trees. 

\medskip \noindent
Now, noticing that $[a_1b_1]_{\mathcal{P}} = \{ a_ip^nb_j \mid i,j \in \{ 1,2,3 \}, \ n \geq 0 \}$ and that a subword $w$ of $a_ip^nb_j$ satisfies $[w]_{\mathcal{P}} \neq \{w\}$ if and only if it contains $a_i$ or $b_j$, we deduce from Proposition \ref{dimSquier} that $\dim S(\mathcal{P},a_1b_1)=2$. $\square$

\begin{remark}
An alternative proof of Corollary \ref{ZbulletZtrees} can be deduced from Example \ref{ZbulletZspecial}, by noticing that $\mathbb{Z} \bullet \mathbb{Z}$ is an undistorted subgroup of $\mathbb{F}_2 \times \mathbb{F}_2$ (it can be proved using the normal form associated to the decomposition of $\mathbb{Z} \bullet \mathbb{Z}$ as an HNN extension given in Example \ref{ZbulletZgraph}).
\end{remark}

\section{Squier complexes as a graph of spaces}

\noindent
For convenience, we begin this section by giving the precise definitions of \textit{graphs of spaces} and \textit{graphs of groups}. For more information, see \cite{MR564422}.

\begin{definition}
A \textit{graph of spaces} is a graph $\Gamma$ such that
\begin{itemize}
	\item each vertex $v \in V(\Gamma)$ is labelled by a space $S_v$,
	\item each edge $e \in E(\Gamma)$ is labelled by a space $S_e$,
	\item for each edge $e=(e^-,e^+)$, there are two $\pi_1$-injective gluing maps
\begin{center}
				$p_e^{\pm} : S_e \to S_{e^{\pm}}$
\end{center}
\end{itemize}
Often, a graph of spaces is identified with its \textit{geometric realization}
\begin{center}
$\displaystyle \left(\bigcup\limits_{v \in V(\Gamma)} S_v \cup \bigcup\limits_{e \in E(\Gamma)} S_e \times [0,1] \right) / \sim$
\end{center}
where $\sim$ identify $S_e \times \{0 \}$ with the image of $p_{e^-}$ and $S_e \times \{ 1 \}$ with the image of $p_{e^+}$.
\end{definition}

\begin{definition}
A \textit{graph of groups} is a graph $\Gamma$ such that
\begin{itemize}
	\item each vertex $v \in V(\Gamma)$ is labelled by a group $G_v$,
	\item each edge $e \in E(\Gamma)$ is labelled by a group $G_e$,
	\item for each edge $e=(e^-,e^+)$, we have two monomorphisms
\begin{center}
				$\varphi_e^{\pm} : G_e \hookrightarrow G_{e^{\pm}}$.
\end{center}
\end{itemize}
Then, we define its \textit{fundamental group} as
\begin{center}
$\left( \underset{v \in V(\Gamma)}{\ast} G_v \ast E(\Gamma) \right) / \left\langle \left\langle E(T), \varphi_{e^-}(g)^{-1}t \varphi_{e^+}(g)t^{-1} \ \text{for all} \ g \in G_e   \right\rangle \right\rangle$,
\end{center}
where $T \subset \Gamma$ is a maximal subtree. According to \cite[Proposition I.5.1.20]{MR1954121}, the group does not depend on the choice of $T$.
\end{definition}

\noindent
Notice that, to any graph of spaces, is naturally associated a graph of groups with the same underlying graph, where the vertex-groups and edge-groups are respectively the fundamental groups of the vertex-spaces and edge-spaces. Then

\begin{thm}
Let $X$ be a connected graph of spaces. The fundamental group of $X$ and the fundamental group of the associated graph of groups coincide.
\end{thm}

\subsection{Decomposition theorem}

\noindent
We fix a semigroup presentation $\mathcal{P}= \langle \Sigma \mid \mathcal{R} \rangle$ and a base word $w \in \Sigma^+$.

\begin{definition}
A hyperplane $J= [a,u \to v,b]$ in $S(\mathcal{P},w)$ is \textit{left} if $D(\mathcal{P},a)= \{1 \}$ but $D(\mathcal{P},au) \neq \{1 \}$.
\end{definition}

\noindent
Notice that $S(\mathcal{P},w)$ contains a left hyperplane whenever the diagram group $D(\mathcal{P},w)$ is not trivial. In this section, we will use Corollary \ref{morphism1} and Lemma \ref{morphism2} intensively, without mentioning them explicitly.

\medskip \noindent
The following two lemmas show that left hyperplanes have good intersection properties. 

\begin{lemma}\label{lefthyperplanesareclean}
A left hyperplane does not self-intersect nor self-osculate, ie., is clean.
\end{lemma}

\noindent
\textbf{Proof.} Let $J= [a,u \to v,b]$ be a left hyperplane. If $J$ self-intersects, we deduce from Lemma \ref{hypautointersectionSquier} and Lemma \ref{hyperplansSquier} that there exists $c \in \Sigma^+$ such that the equalities $a=auc$ and $b=cub$ hold modulo $\mathcal{P}$. Therefore,
\begin{center}
$\{ 1 \} \neq D(\mathcal{P},au) \hookrightarrow D(\mathcal{P},auc) \simeq D(\mathcal{P},a) = \{1 \}$,
\end{center}
a contradiction. Then, if $J$ self-osculates, we deduce from Lemma \ref{hypautotangenceSquier} and Lemma \ref{hyperplansSquier} that there exist $n \geq 1$ and $h,k \in \Sigma^+$ such that the equality $u=(kh)^nk$ holds in $\Sigma^+$ and the equalities $a=akh$ and $b=hkb$ hold modulo $\mathcal{P}$. Similarly, we get
\begin{center}
$\{ 1 \} \neq D(\mathcal{P},au) = D(\mathcal{P},a(kh)^nk) \hookrightarrow D(\mathcal{P},a(kh)^{n+1}) \simeq D(\mathcal{P},a) = \{1 \}$,
\end{center}
a contradiction. $\square$

\begin{lemma}
Two left hyperplanes do not intersect.
\end{lemma}

\noindent
\textbf{Proof.} Let $J_1,J_2$ be two left hyperplanes. If they intersect, there exists a square
\begin{center}
$(a,u \to v, b, p \to q,c)$
\end{center}
such that, say, $J_1=[a,u \to v, bpc]$ and $J_2= [aub,p \to q,c]$. We deduce that
\begin{center}
$\{1 \} \neq D(\mathcal{P},au) \hookrightarrow D(\mathcal{P},aub) = \{1 \}$,
\end{center}
because $J_1$ and $J_2$ are left, whence a contradiction. $\square$

\medskip \noindent
If $J= [a,u \to v,b]$ is a left hyperplane, let $p_J,s_J \in \Sigma^+$ be two words and $\ell_J \in \Sigma$ a letter satisfying: $u=p_J \ell_J s_J$ in $\Sigma^+$, $D(\mathcal{P},ap_J)= \{1 \}$ and $D(\mathcal{P},ap_J \ell_J) \neq \{1 \}$. Notice that $p_J$ is just the maximal prefix of $u$ satisfying $D(\mathcal{P},ap_J)= \{1 \}$, so that $p_J,s_J, \ell_J$ are uniquely determined. We define similarly $v=q_Jm_Jr_J$, where $q_J,r_J \in \Sigma^+$ and $m_J \in \Sigma$, so that $D(\mathcal{P},q_J)= \{1 \}$ and $D(\mathcal{P},q_Jm_J) \neq \{1 \}$.

\medskip \noindent
This notation is motivated by the following technical lemma:

\begin{lemma}\label{technicallemma}
Let $J=[a,u \to v,b]$ be a left hyperplane. Let $x,y,p,q \in \Sigma^+$ be such that $aub=xpy$ in $\Sigma^+$ and $p=q \in \mathcal{R}$. Then, $H=[x,p \to q,y]$ is not a left hyperplane if and only if $p$, as a subword of $aub$, is included into $ap_J$ or $s_Jb$. 
\end{lemma}

\noindent
\textbf{Proof.} If $p$ is included into $ap_J$, then 
\begin{center}
$D(\mathcal{P},xp) \hookrightarrow D(\mathcal{P},ap_J)= \{1 \}$,
\end{center}
so $H$ is not a left hyperplane. If $p$ is included into $s_Jb$, then $ap_J\ell_J$ is included into $x$, hence
\begin{center}
$\{1 \} \neq D(\mathcal{P},ap_J \ell_J) \hookrightarrow D(\mathcal{P},x)$,
\end{center}
so $H$ is not a left hyperplane. 

\medskip \noindent
Conversely, suppose that $p$ is included neither into $ap_J$ nor into $s_Jb$. Then $p$ contains the letter $\ell_J$, viewed as a subword of $u$; as a consequence, $x$ is included into $ap_J$, hence
\begin{center}
$D(\mathcal{P},x) \hookrightarrow D(\mathcal{P},ap_J)= \{1 \}$
\end{center}
and
\begin{center}
$\{1 \} \neq D(\mathcal{P},ap_J \ell_J ) \hookrightarrow D(\mathcal{P},xp)$.
\end{center}
Therefore, $H$ is a left hyperplane. $\square$

\medskip \noindent
Now, we are ready to state and prove the main theorem of this section. Roughly speaking, we find a graph of spaces by cutting $S(\mathcal{P},w)$ along its left hyperplanes.

\begin{thm}\label{graphofspaces}
Let $\mathfrak{J}$ denote the set of left hyperplanes of $S(\mathcal{P},w)$. Let $\mathcal{G}(\mathcal{P},w)$ be the graph of spaces defined by:
\begin{itemize}
	\item the set of vertex-spaces is
\begin{center}
$\{ S(\mathcal{P},ap_J) \ell_J S(\mathcal{P},s_Jb), \ S(\mathcal{P},aq_J) m_J S(\mathcal{P},r_Jb) \mid J= [a,u \to v,b] \in \mathfrak{J} \}$,
\end{center}
	\item to each left hyperplane $[a, u \to v,b] \in \mathfrak{J}$ is associated the edge-space $S(\mathcal{P},a) \times S(\mathcal{P},b)$,
	\item the edge-maps are the canonical maps
\begin{center}
$S(\mathcal{P},a) \times S(\mathcal{P},b) \to S(\mathcal{P},a)u S(\mathcal{P},b)$ and $S(\mathcal{P},a) \times S(\mathcal{P},b) \to S(\mathcal{P},a)v S(\mathcal{P},b)$.
\end{center}
\end{itemize}
Then $\mathcal{G}(\mathcal{P},w)$ defines a decomposition of $S(\mathcal{P},w)$ as a graph of spaces.
\end{thm}

\begin{remark}\label{remark}
We emphasize that, if two distinct left hyperplanes define two identical vertex-spaces, then these spaces define only one vertex in the graph of spaces. However, if two distinct left hyperplanes define two identical edge-spaces, then these spaces define two edges in the graph of spaces.
\end{remark}

\noindent
\textbf{Proof of Theorem \ref{graphofspaces}.} Let $\overline{S}(\mathcal{P},w)$ denote the subcomplex
\begin{center}
$S(\mathcal{P},w) \backslash \bigcup\limits_{J \in \mathfrak{J}} \left( N(J) \backslash \partial J \right)$.
\end{center}
Then, since left hyperplanes are clean according to Lemma \ref{lefthyperplanesareclean}, $S(\mathcal{P},w)$ is constructed from the connected components of $\overline{S}(\mathcal{P},w)$ by taking a copy of $N(J) \simeq J \times [0,1]$ for each left hyperplane $J \in \mathfrak{J}$ and gluing $J \times \{0 \}$ and $J \times \{1 \}$ along $\partial_- J$ and $\partial_+ J$ respectively via the natural isomotries $J \to \partial_- J$ and $J \to \partial_+ J$ given by Lemma \ref{partialJ}. In particular, notice that these gluings are $\pi_1$-injective, so that $S(\mathcal{P},w)$ may be decomposed as the graph of spaces defined by:
\begin{itemize}
	\item the vertices are the connected components of $\overline{S}(\mathcal{P},w)$,
	\item to each left hyperplane $J$ is associated an edge linking the two (not necessarily distinct) connected components adjacent to $J$,
	\item if $J= [a,u \to v,b]$ is a left hyperplane, the gluing maps $J \times \{0 \} \to \partial_- J$ and $J \times \{ 1 \} \to  \partial_+J$ are given by 
\begin{center}
	$J \simeq S(\mathcal{P},a) \times S(\mathcal{P},b) \to S(\mathcal{P},a)uS( \mathcal{P},b)= \partial_- J$
\end{center}
				and
\begin{center}
	$J \simeq S(\mathcal{P},a) \times S(\mathcal{P},b) \to S(\mathcal{P},a)vS( \mathcal{P},b)= \partial_+ J$
\end{center}
				respectively, following Lemma \ref{partialJ}.	
\end{itemize}
Therefore, to conclude the proof, it is sufficient to prove that, if $J= [a,u \to v,b]$ is a left hyperplane and $C_u$ (resp. $C_v$) the connected component of $\overline{S}(\mathcal{P},w)$ containing $aub$ (resp. $avb$), then $C_u= S(\mathcal{P},ap_J)\ell_J S(\mathcal{P},s_Jb)$ (resp. $C_v= S(\mathcal{P},aq_J)m_J S(\mathcal{P},r_Jb)$). In fact, by symmetry, we only have to prove the claim for $C_u$. 

\medskip \noindent
Let $e$ be an oriented edge of $C_u$. We prove by induction on the combinatorial distance $d$ between $aub \in C_u$ and the initial point $e^-$ of $e$ that $e$ belongs to $S(\mathcal{P},ap_J)\ell_J S(\mathcal{P},s_Jb)$. If $d=0$, then $e^-=aub$, so that $e=(x,p \to q,y)$ for some words $x,y,p,q \in \Sigma^+$ satisfying $xpy=aub$ in $\Sigma^+$. Because $e$ belongs to $C_u$, by definition $[x,p \to q,y]$ cannot be a left hyperplane, so, according to Lemma \ref{technicallemma}, $p$ is included into $ap_J$ or $s_Jb$. It follows that $e$ is an edge of $S(\mathcal{P},ap_J)\ell_J S(\mathcal{P},s_Jb)$. Now, let $d \geq 1$. If we consider a path of edges linking $aub$ to $e$, by the induction hypothesis, we know that the penultimate edge belongs to $S(\mathcal{P},ap_J)\ell_J S(\mathcal{P},s_Jb)$, so $e^-= \alpha \ell_J \beta$ for some words $\alpha, \beta \in \Sigma^+$ satisfying $\alpha=ap_J$ and $\beta=s_Jb$ modulo $\mathcal{P}$. Let $x,y,p,q \in \Sigma^+$ be such that $e=(x,p \to q,y)$; in particular, $xpy=e^-=\alpha \ell_J \beta$ in $\Sigma^+$. If $p$, considered as a subword of $\alpha \ell_J \beta$, contains the letter $\ell_J$, then
\begin{center}
$D(\mathcal{P},x) \hookrightarrow D(\mathcal{P}, \alpha) \simeq D(\mathcal{P},ap_J)= \{1 \}$
\end{center}
and
\begin{center}
$\{1 \} \neq D(\mathcal{P},ap_J\ell_J) \simeq D(\mathcal{P}, \alpha \ell_J) \hookrightarrow D(\mathcal{P},xp)$,
\end{center}
so $e$ is dual to a left hyperplane: a contradiction with $e \subset C_u$. Therefore, $p$ is included into $\alpha$ or $\beta$, so that $e=(x,p \to q,y)$ belongs to $S(\mathcal{P},ap_J)\ell_J S(\mathcal{P},s_Jb)$.

\medskip \noindent
Thus, we have proved that $C_u \subset S(\mathcal{P},ap_J)\ell_J S(\mathcal{P},s_Jb)$. Conversely, we deduce from Lemma \ref{technicallemma} that no edge of $S(\mathcal{P},ap_J)\ell_J S(\mathcal{P},s_Jb)$ is dual to a left hyperplane, so that the inclusion $S(\mathcal{P},ap_J)\ell_J S(\mathcal{P},s_Jb) \subset C_u$ holds. $\square$

\begin{ex}\label{ZbulletZgraph}
Let $\mathcal{P}= \left\langle \begin{array}{l} a_1,a_2,a_3 \\ b_1,b_2,b_3 \end{array}, p \left| \begin{array}{l} a_1=a_2, a_2=a_3, a_3=a_1 \\ b_1=b_2, b_2=b_3, b_3=b_1 \end{array}, \begin{array}{l} a_1=a_1p \\ b_1=pb_1 \end{array} \right. \right\rangle$.

\medskip \noindent
Then $S(\mathcal{P},a_1b_1)$ contains four left hyperplanes: $[1,a_1 \to a_2, b_1]$, $[1,a_2 \to a_3,b_1]$, $[1,a_3 \to a_1,b_1]$ and $[1,a_1 \to a_1p,b_1]$. Thus, the vertex-spaces of our graph of spaces will be $a_1S(\mathcal{P},b_1)$, $a_2S(\mathcal{P},b_1)$ and $a_3S(\mathcal{P},b_1)$. The graph of spaces given by Theorem \ref{graphofspaces} is:
\begin{center}
\includegraphics[scale=0.5]{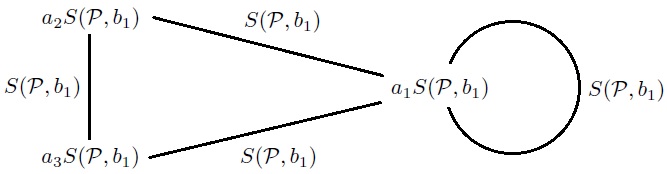}
\end{center}
\noindent
The maps associated to the loop are induced by $w_0 \mapsto a_1w_0$ and $w_0 \mapsto a_1pw_0$. On the other hand, the Squier complex $S(\mathcal{P},b_1)$ is easy to draw:
\begin{center}
\includegraphics[scale=0.5]{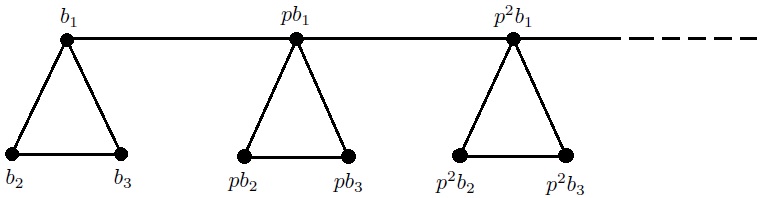}
\end{center}
In particular, its fundamental group is isomorphic to $\mathbb{F}_{\infty}= \langle x_1,x_2, \ldots \mid \ \rangle$. Thus, $D(\mathcal{P},a_1b_1)$ may be decomposed as the following graph of groups:
\begin{center}
\includegraphics[scale=0.6]{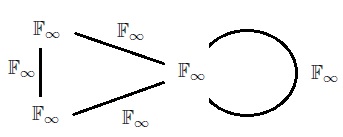}
\end{center}
\noindent
where the maps associated to the three edges on the left are identities, and where the two maps associated to the loop are the identity and the morphism induced by $x_i \mapsto x_{i+1}$. This graph of groups may be simplified into:
\begin{center}
\includegraphics[scale=0.6]{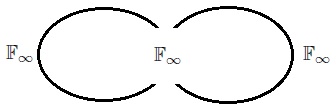}
\end{center}
\noindent
where the maps associated to the loop on the left are identities, and where the two maps associated to the loop on the right are the identity and the morphism induced by $x_i \mapsto x_{i+1}$. Thus, a presentation of $D(\mathcal{P},a_1b_1)$ is:
\begin{center}
$\langle x_1,x_2, \ldots, t,h \mid tx_it^{-1}=x_i, hx_ih^{-1}=x_{i+1} \ (i \geq 1 ) \rangle$.
\end{center}
Noticing that $x_{i+1}=x_1^{h^i}$, we simplify the presentation above into
\begin{center}
$\mathbb{Z} \bullet \mathbb{Z} = \langle a,t,h \mid [t,a^{h^i}]=1 \ (i \geq 0) \rangle$.
\end{center}
\end{ex}

\begin{ex}
Let $\mathcal{P}= \langle a,b,c \mid ab=ba, ac=ca, bc=cb \rangle$. For convenience, let $U(\ell,m,n)$ denote the diagram group $D(\mathcal{P},a^{\ell}b^mc^n)$. We gave in \cite[Example 5.3]{arXiv:1505.02053} a criterion to determine whether or not a given diagram group is free, and we noticed that $U(\ell,m,n)$ is free whenever $\ell$, $m$ or $n$ is $1$. However, this method did not give any information on the rank of the free group. Using Theorem \ref{graphofspaces}, we are now able to prove that $U(1,m,n)$ is a free group of rank $mn$. (This result is also proved in \cite[Example 10.2]{MR1396957}.)

\medskip \noindent
The left hyperplanes of $S(\mathcal{P},ab^mc^n)$ are:
\begin{itemize}
	\item $J_{k \ell}^b = [b^kc^{\ell}, ab \to ba, b^{m-k-1}c^{n-\ell}]$ with $\left\{ \begin{array}{l} 0 \leq k \leq m-1 \\ 1 \leq \ell \leq n \end{array} \right.$,
	\item $J_{k \ell}^c = [b^kc^{\ell}, ac \to ca, b^{m-k}c^{n- \ell-1}]$ with $\left\{ \begin{array}{l} 1 \leq k \leq m \\ 0 \leq \ell \leq n-1 \end{array} \right.$,
	\item $H^b_k= [ ab^k,bc \to cb, b^{m-k-1}c^{n-1}]$ with $0 \leq k \leq m-1$,
	\item $H^c_k= [ac^k,bc \to cb, b^{m-1}c^{n-k-1}]$ with $0 \leq k \leq n-1$.
\end{itemize}
It is worth noticing that $H^b_0=H^c_0$ is the only hyperplane appearing more than once in the list above so that there are $2mn+m+n-1$ left hyperplanes. Then, in our graph of spaces,
\begin{itemize}
	\item $J_{0 \ell}^b$ links $S(\mathcal{P},ac^{\ell})bS(\mathcal{P},b^{m-1}c^{n-\ell})$ to $S(\mathcal{P},bc^{\ell})aS(\mathcal{P}, b^{m-1}c^{n- \ell})$,
	\item if $k \geq 1$, $J_{k \ell}^b$ links $S(\mathcal{P},b^kc^{\ell})aS(\mathcal{P},b^{m-k}c^{n- \ell})$ to $S(\mathcal{P},b^{k+1} c^{\ell})aS(\mathcal{P}, b^{m-k-1}c^{n-\ell})$,
	\item $J_{0 \ell}^c$ links $S(\mathcal{P},ab^{\ell})cS(\mathcal{P},b^{m-\ell}c^{n-1})$ to $S(\mathcal{P},b^{\ell} c)aS(\mathcal{P}, b^{m-\ell}c^{n-1})$,
	\item if $k \geq 1$, $J_{k \ell}^c$ links $S(\mathcal{P},b^{\ell}c^{k})aS(\mathcal{P},b^{m-\ell}c^{n-k})$ to $S(\mathcal{P},b^{\ell} c^{k+1})aS(\mathcal{P}, b^{m-\ell}c^{n-k-1})$,
	\item $H_0^b$ links $S(\mathcal{P},ab)cS(\mathcal{P},b^{m-1}c^{n-1})$ to $S(\mathcal{P},ac)bS(\mathcal{P},b^{m-1}c^{n-1})$,
	\item if $k \geq 1$, $H_k^b$ links $S(\mathcal{P},ab^{k+1})cS(\mathcal{P},b^{m-k-1}c^{n-1})$ to $S(\mathcal{P},ab^k)cS (\mathcal{P}, b^{m-k} c^{n-1})$,
	\item if $k \geq 1$, $H_k^c$ links $S(\mathcal{P},ac^{k+1})bS(\mathcal{P},b^{m-1}c^{n-k-1})$ to $S(\mathcal{P},ac^k)bS (\mathcal{P}, b^{m-1} c^{n-k})$,
\end{itemize}
In the list above, $mn+m+n$ vertex-spaces appear. 

\medskip \noindent
Now, since $S(\mathcal{P},w)$ is simply connected when the word $w$ has at most two different letters, we deduce that all the vertex-spaces and edge-spaces are simply connected, so that the vertex-groups and edge-groups in the associated graph of groups are trivial. Thus, $U(1,m,n)$ is the fundamental group of a simplicial graph with $mn+m+n$ vertices and $2mn+m+n-1$ edges: it is a free group of rank
\begin{center}
$(2mn+m+n-1)- (mn+m+n-1) = mn$.
\end{center}
\end{ex}

\begin{remark}
Similarly, \textit{right hyperplanes} may be defined: a hyperplane $[a,u \to v,b]$ is \textit{right} whenever $D(\mathcal{P},b)=\{1 \}$ but $D(\mathcal{P},ub) \neq \{1 \}$. Then Theorem \ref{graphofspaces} has an equivalent statement for right hyperplanes. For example, the decomposition of the Squier complex $S(\mathcal{P}_1,x)$, with the semigroup presentation
\begin{center}
$\mathcal{P}_1 = \langle x,a,b,c \mid x=ax, \ a=b, \ b=c, \ c=a \rangle$,
\end{center}
is more efficient with respect to right hyperplanes: it allows to prove that the diagram group $D(\mathcal{P}_1,x)$ is a free abelian group of infinite (countable) rank. Compare with the example given in Remark \ref{Zinfty}.
\end{remark}

\subsection{Application: right-angled Artin groups and interval graphs}

\noindent
This section is dedicated to the proof of Theorem \ref{mainintervalgraphRAAG}.

\medskip \noindent
Let $\Gamma$ be a finite interval graph. Since it is finite, we may suppose without loss of generality that $\Gamma$ is associated to a collection $\mathcal{C}$ of intervals on $\{1, \ldots,n \}$, for some $n \geq 1$. For convenience, if $I = ( i_1,\ldots, i_r )$, we will note $x_I = x_{i_1} \cdots x_{i_r}$. Then, to the collection $\mathcal{C}$, we associate the semigroup presentation
\begin{center}
$\mathcal{P}(\mathcal{C})= \langle x_1,\ldots,x_n,a_I,b_I,c_I \ (I \in \mathcal{C}) \mid x_I=a_I, \ a_I=b_I, \ b_I=c_I, \ c_I=a_I \ (I \in \mathcal{C}) \rangle$.
\end{center}
The main result of this section is

\begin{thm}\label{principal'}
The diagram group $D(\mathcal{P}(\mathcal{C}), x_1 \cdots x_n)$ is isomorphic to the right-angled Artin group $A(\overline{\Gamma})$.
\end{thm}

\noindent
\textbf{Proof.} For convenience, let $\mathcal{P}= \mathcal{P}(\mathcal{C})$. For all $I \in \mathcal{C}$, let $\Delta_I \in D(\mathcal{P},x_1 \cdots x_n)$ be the following spherical diagram:
\begin{center}
\includegraphics[scale=0.6]{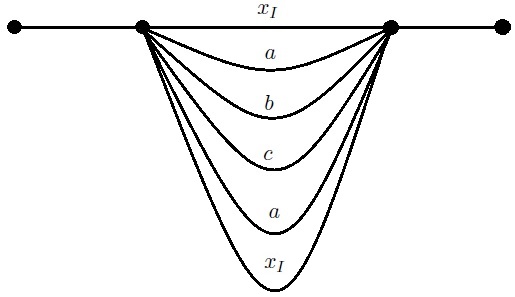}
\end{center}
Notice that, if $I$ and $J$ are disjoint, then $\Delta_I$ and $\Delta_J$ commute. Thus, there is a natural morphism 
\begin{center}
$\Phi : \left\{ \begin{array}{ccc} A( \overline{\Gamma} ) & \to & D(\mathcal{P}(\mathcal{C}),x_1 \cdots x_n)  \\  I & \mapsto & \Delta_I \end{array} \right.$.
\end{center}
We want to prove by induction on the number of vertices of $\Gamma$ that $\Phi$ is an isomorphism. If $\Gamma$ has no vertex, ie., if $\mathcal{C}$ is empty, then the two groups are trivial and there is nothing to prove.

\medskip \noindent
From now on, suppose that $\Gamma$ has at least one vertex, ie., $\mathcal{C}$ contains at least one interval. If $I_1, I_2 \in \mathcal{C}$ are two disjoint intervals and if $I_1$ is at the left of $I_2$, we will note $I_1 \prec I_2$; furthermore, an interval $I \in \mathcal{C}$ will be \textit{left} if $I$ is minimal in $\mathcal{C}$ with respect to $\prec$. Finally, for all $I \in \mathcal{C}$, let $g_I,d_I \in \Sigma^+$ be the words satisfying the equality $g_Ix_Id_I=x_1 \cdots x_n$ in $\Sigma$; notice that $I$ is a left interval if and only if $[g_I,x_I \to a_I,d_I]$ is a left hyperplane.

\medskip \noindent
We want to decompose the Squier complex $S(\mathcal{P},x_1 \cdots x_n)$ as a graph of spaces thanks to Theorem \ref{graphofspaces}. The left hyperplanes are 
\begin{itemize}
	\item $X_I=[g_I,x_I \to a_I,d_I]$, 
	\item $A_I=[g_I,a_I \to b_I,d_I]$, 
	\item $B_I=[g_I, b_I \to c_I,d_I]$,
	\item $C_I=[g_I,c_I \to a_I, d_I]$,
\end{itemize}
for all left interval $I$. Let $x_I=p_I \ell_Is_I$ be the decomposition used in Theorem \ref{graphofspaces}, ie., $p_I,s_I \in \Sigma^+$ and $ \ell_I \in \Sigma$ satisfy $D(\mathcal{P},g_Ip_I)= \{1 \}$ and $D(\mathcal{P},g_Ip_I \ell_I) \neq \{ 1 \}$. Then, in our graph of spaces,
\begin{itemize}
	\item $X_I$ links $S(\mathcal{P},g_Ip_I) \ell_I S(\mathcal{P},s_Id_I)$ and $S(\mathcal{P},g_I)a_I S(\mathcal{P},d_I)$,
	\item $A_I$ links $S(\mathcal{P},g_I)a_IS(\mathcal{P},d_I)$ and $S(\mathcal{P},g_I)b_IS(\mathcal{P},d_I)$,
	\item $B_I$ links $S(\mathcal{P},g_I)b_IS(\mathcal{P},d_I)$ and $S(\mathcal{P},g_I)c_IS(\mathcal{P},d_I)$,
	\item $C_I$ links $S(\mathcal{P},g_I)c_IS(\mathcal{P},d_I)$ and $S(\mathcal{P},g_I)a_IS(\mathcal{P},d_I)$,
\end{itemize}
Notice that $S(\mathcal{P},g_Ip_I)=\{g_Ip_I\}$ and $S(\mathcal{P},g_I)=\{g_I\}$. Indeed, if it was not the case, there would exist a $x_J$ included into $g_Ip_I$ or $g_I$ and we would deduce that $D(\mathcal{P},g_Ip_I) \neq \{1 \}$ or $D(\mathcal{P},g_I) \neq \{1 \}$, which is in contradiction with the fact that $[g_I,x_I,d_I]$ is a left hyperplane. Therefore,
\begin{itemize}
	\item $X_I$ links $g_Ip_I \ell_I S(\mathcal{P},s_Id_I)$ and $g_Ia_I S(\mathcal{P},d_I)$,
	\item $A_I$ links $g_Ia_IS(\mathcal{P},d_I)$ and $g_Ib_IS(\mathcal{P},d_I)$,
	\item $B_I$ links $g_Ib_IS(\mathcal{P},d_I)$ and $g_Ic_IS(\mathcal{P},d_I)$,
	\item $C_I$ links $g_Ic_IS(\mathcal{P},d_I)$ and $g_Ia_IS(\mathcal{P},d_I)$,
\end{itemize}
Thus, our vertex-spaces are $g_Ip_I \ell_I S(\mathcal{P},s_Id_I)$, $g_Ia_IS(\mathcal{P},d_I)$, $g_Ib_IS(\mathcal{P},d_I)$ and $g_Ic_IS(\mathcal{P},d_I)$ when $I$ runs over the set of left intervals of $\mathcal{C}$. However, according to Remark \ref{remark}, we have to compare these spaces. The only non-trivial question is to determine whether or not $g_Ip_I \ell_I S(\mathcal{P},s_Id_I)$ and $g_Jp_J \ell_J S(\mathcal{P},s_Jd_J)$ are different. 

\medskip \noindent
We want to prove that the equalities $g_Ip_I \ell_I=g_J p_J \ell_J$ and $s_Id_I=s_Jd_J$ hold in $\Sigma^+$ for any left hyperplanes $I,J$. As a consequence, we will deduce that the spaces $g_Ip_I \ell_I S(\mathcal{P},s_Id_I)$ and $g_Jp_J \ell_J S(\mathcal{P},s_Jd_J)$ define only one vertex in our graph of spaces. Notice first that, because 
\begin{center}
$g_Ip_I\ell_Is_Id_I= x_1 \cdots x_n= g_Jp_J \ell_Js_J d_J$
\end{center}
in $\Sigma^+$, it is sufficient to prove that $g_Ip_I \ell_I=g_J p_J \ell_J$: the second equality then follows.

\medskip \noindent
If $g_Ip_I \ell_I \neq g_J p_J \ell_J$ then either $g_Ip_I \ell_I$ is a proper prefix of $g_Jp_J \ell_J$ or $g_Jp_J \ell_J$ is a proper prefix of $g_Ip_I \ell_I$. In the first case, we would have
\begin{center}
$\{1 \} \neq D(\mathcal{P},g_Ip_I \ell_I) \hookrightarrow D(\mathcal{P},g_Jp_J)= \{1 \}$,
\end{center}
and similarly, in the second case we would have
\begin{center}
$\{1 \} \neq D(\mathcal{P},g_Jp_J \ell_J) \hookrightarrow D(\mathcal{P},g_Ip_I)= \{1 \}$.
\end{center}
Therefore, we conclude that $g_Ip_I \ell_I = g_J p_J \ell_J$.

\medskip \noindent
In particular, we may note $g=g_Ip_I \ell_I$ and $d=s_Id_I$ so that $g$ and $d$ do not depend on $I$. Let $I_1,\ldots, I_r$ be the left intervals of $\mathcal{C}$. Now, our graph of spaces is
\begin{center}
\includegraphics[scale=0.59]{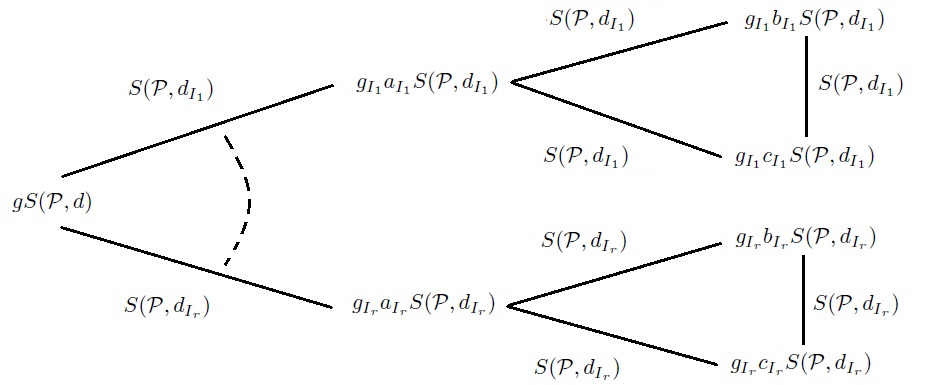}
\end{center}
Notice that, using the morphism of Lemma \ref{iso}, $\Delta_{I_k}$ defines a loop starting from $x_1 \cdots x_n \in gS( \mathcal{P},d)$ and passing through the $k$-th loop of our graph of spaces: the $\Delta_{I_k}$ will define the stable letters of the HNN extensions.

\medskip \noindent
Noticing that the edge-maps of the graph of spaces are all identities, we deduce that $D(\mathcal{P},x_1 \cdots x_n)$ is an HNN extension of $\pi_1(g S(\mathcal{P},d),x_1 \cdots x_n)$ over the subgroups 
\begin{center}
$\pi_1(g_{I_k}x_{I_k}S(\mathcal{P},d_{I_k}),x_1 \cdots x_n)$ 
\end{center}
with stable letter $\Delta_{I_k}$. 

\medskip \noindent
Because $d$ is a subword of $x_1 \cdots x_n$, there exists an interval $J \subset \{ 1, \ldots, n \}$ such that $d=x_J$. Now, notice that a diagram $\Delta$ satisfying $\mathrm{top}(\Delta)=d$ can contain a cell corresponding to a relation of $\mathcal{P}$ of the form $x_I \to a_I$  if and only if $J \subset I$. Therefore, if we introduce the following semigroup presentation,
\begin{center}
$\mathcal{P}_d= \left\langle \begin{array}{l} a_I,b_I,c_I, x_k, \\ (k \in J, I \in \mathcal{C} \ \text{and} \ I \subset J) \end{array} \left| \begin{array}{l} a_I=b_I, \ b_I=c_I, \\ c_I=a_I, \ x_I=a_I, \end{array} \ (I \subset J \ \text{and} \ I \in \mathcal{C}) \right. \right\rangle$,
\end{center}
then $gS(\mathcal{P},d)=S(\mathcal{P}_d,gd)$. By our induction hypothesis, if $\Gamma_0$ is the subgraph of $\Gamma$ generated by the vertices corresponding to the intervals of $\mathcal{C}$ which are not left, then the map 
\begin{center}
$\left\{ \begin{array}{ccc} A( \overline{\Gamma_0} ) & \to & D(\mathcal{P}_d,gd) \\ I & \mapsto & \Delta_I \end{array} \right.$
\end{center}
defines an isomorphism. Similarly, the fundamental group of 
\begin{center}
$g_{I_k}x_{I_k} S(\mathcal{P},d_{I_k})=g_{I_k}x_{I_k} S(\mathcal{P}_d,d_{I_k})$ 
\end{center}
coincides with the subgroup of the fundamental group of $S(\mathcal{P}_d,gd)$ generated by $\{ \Delta_J \mid J \cap I_k = \emptyset \}$.

\medskip \noindent
Consequently, the fundamental group of the Squier complex $S(\mathcal{P},x_1 \cdots x_n)$, ie., the diagram group $D(\mathcal{P},x_1 \cdots x_n)$, is an HNN extension over the subgroup $\langle \Delta_I \mid I \ \text{is not left} \rangle$, which is isomorphic to the right-angled Artin group $A( \overline{\Gamma_0})$, where the stable letters are $\{ \Delta_I \mid I \ \text{is left} \}= \{ \Delta_{I_1}, \ldots, \Delta_{I_k} \}$ and where each $\Delta_{I_k}$ has to commute with $\{ \Delta_I \mid I \cap I_k = \emptyset \}$. This description exactly means that the morphism $\Phi$, from $A(\Gamma)$ to $D(\mathcal{P},x_1 \cdots x_n)$, is an isomorphism. $\square$

\begin{remark}
The semigroup presentation $\mathcal{P}(\mathcal{C})$ used above is \textit{complete}, so that it is possible to apply the algorithm \cite[Theorem 9.8]{MR1396957} to find a presentation of the diagram group $D(\mathcal{P}(\mathcal{C}),x_1 \cdots x_n)$. The presentation we find is exactly the canonical presentation of $A(\overline{\Gamma})$, ie.,
\begin{center}
$\langle A_I \ (I \in \mathcal{C}) \mid [A_I,A_J]=1 \ (I \cap J = \emptyset) \rangle$.
\end{center}
Thus, this gives an alternative proof of Theorem~\ref{principal'}. 
\end{remark}

\noindent
Of course, a natural question follows from Theorem \ref{principal'}: when is a finite graph the complement of an interval graph? A simple criterion is given in \cite[Theorem 3.5]{MR1672910}.

\begin{definition}
Let $\Gamma$ be a graph and let $V(\Gamma)$ (resp. $E(\Gamma)$) denote the set of vertices (resp. edges) of $\Gamma$. The graph $\Gamma$ is \textit{transitively orientable} if it admits an orientation satisfying the following property: for any vertices $x,y,z \in V(\Gamma)$, $(x,y) \in E(\Gamma)$ and $(y,z) \in E(\Gamma)$ implies $(x,z) \in E(\Gamma)$.
\end{definition}

\begin{thm}
A graph is the complement of an interval graph if and only if it does not contain $\overline{C_4}$ as an induced subgraph and if it is transitively orientable.
\end{thm}

\noindent
Recall that $\Gamma \subset \Lambda$ is an \textit{induced subgraph} of $\Lambda$ if any vertices $x,y \in \Gamma$ are linked by an edge in $\Lambda$ if and only if they are linked by an edge in $\Gamma$.

\begin{ex}\label{AP26}
As a consequence, we may deduce that the following graphs are the complements of interval graphs, so that the associated right-angled Artin groups are diagral groups according to Theorem \ref{principal'}. (The blue arrows induce a transitive orientation on the graphs.)
\begin{center}
\includegraphics[scale=0.5]{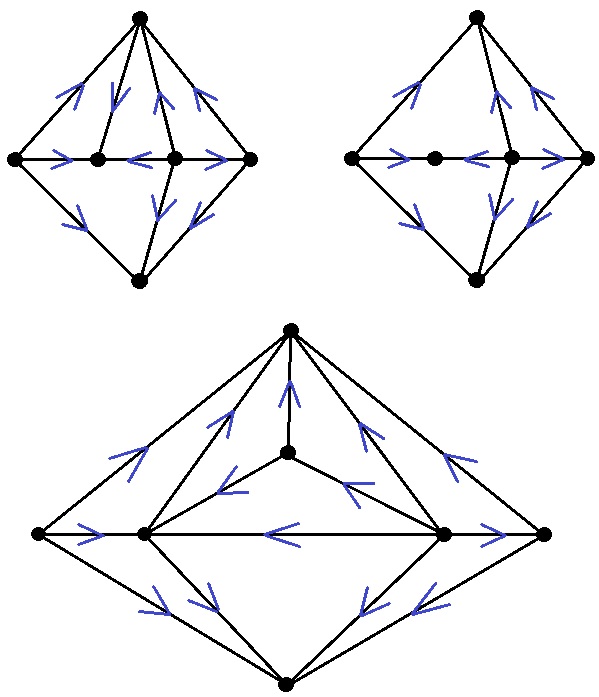}
\end{center}
The first of the three graphs above is denoted by $P_2(6)$ in \cite[Section 7]{MR2422070}; it is proved in \cite[Section 8]{MR2422070} that $A(P_2(6))$ embeds into a diagram group. Thus, Theorem~\ref{principal'} gives a stronger conclusion: it is a diagram group itself.
\end{ex}

\section{Some open questions}

\noindent
In Section 3.3, we constructed a morphism $\Phi= \Phi(\mathcal{P},w) : D(\mathcal{P},w) \to A(\mathcal{P},w)$ for every diagram group $D(\mathcal{P},w)$, where $A(\mathcal{P},w)$ is some right-angled Artin group (see Definition \ref{Phi}), and we proved that $\Phi$ is injective whenever the Squier complex $S(\mathcal{P},w)$ is special (Proposition \ref{Phispecial}). It is worth noticing that the specialness of $S(\mathcal{P},w)$ depends on the semigroup presentation $\mathcal{P}$ and is not an algebraic invariant of the diagram group $D(\mathcal{P},w)$: if
\begin{center}
$\mathcal{P}_1 = \langle a,b,p,q \mid a=ap,b=pb,p=q \rangle$ and $\mathcal{P}_2 = \langle a,b,c \mid a=b,b=c,c=a \rangle$,
\end{center}
then $D(\mathcal{P}_1,ab) \simeq \mathbb{Z} \simeq D(\mathcal{P}_2,a)$, whereas $S(\mathcal{P}_1,ab)$ is not a special cube complex but $S(\mathcal{P}_2,a)$ is. However, both $\Phi(\mathcal{P}_1,ab)$ and $\Phi(\mathcal{P}_2,a)$ are injective. In fact, we suspect that the injectivity of $\Phi$ depends only on the isomorphic class of the diagram group, ie., is an algebraic invariant; and we hope to deduce a proof of the following result:

\begin{conj}\label{embedRAAG}
A diagram group is embeddable into a right-angled Artin group if and only if it does not contain $\mathbb{Z} \wr \mathbb{Z}$. 
\end{conj}

\noindent
This conjecture is motivated by Theorem \ref{Squierspecial}, and by \cite[Theorem 24]{MR1725439} which proves that the diagram group $D(\mathcal{P},w)$ contains $\mathbb{Z} \wr \mathbb{Z}$ if and only if there exist some words $a,b,p \in \Sigma^+$ satisfying $w=ab$, $a=ap$ and $b=pb$ modulo $\mathcal{P}$ with $D(\mathcal{P},p) \neq \{1 \}$. Notice also that a corollary would be that any simple diagram group contains $\mathbb{Z} \wr \mathbb{Z}$: this result is indeed true, and can be deduced from \cite[Theorem 7.2]{MR2193191} and \cite[Corollary 22]{MR1725439}.

\medskip \noindent
A positive answer to Conjecture \ref{embedRAAG} would have several nice corollaries on diagram groups; in particular, it would give a simple criterion of linearity (and a fortiori of residual finiteness). However, such a criterion cannot be necessary, since $\mathbb{Z} \wr \mathbb{Z}$ is linear. However, we ask the following (possibly naive) question:

\begin{question}
Is a diagram group not containing Thompson's group $F$ linear or residually finite?
\end{question}

\noindent
According to \cite[Theorem 4.2]{MR2377497}, a group is embeddable into a right-angled Artin group if and only if it is the fundamental group of a special cube complex. Another interesting question is to know when a diagram group is the fundamental group of a compact special cube complex, in order to apply \cite[Corollary 7.9]{MR2377497} and thus deduce separability of convex-cocompact subgroups. Of course, a necessary condition is to be finitely-presented. We ask whether this condition is sufficient:

\begin{question}
Let $G$ be a diagram group embeddable into a right-angled Artin group. If $G$ is finitely-presented, is it the fundamental group of a compact special cube complex?
\end{question}

\noindent
A solution could be to find an interesting hierarchy of such diagram groups by studying their actions on their associated Farley complexes. 

\medskip \noindent
Conversely, we wonder which right-angled Artin groups are embeddable into a diagram group. Corollary \ref{AC5} proves that the best candidates we had as right-angled Artin groups not embeddable into a diagram group turn out to embed in such a group.

\begin{question}
Does any (finitely generated) right-angled Artin group embed into a diagram group?
\end{question}

\noindent
A positive answer would follow from the fact that any (finitely generated) right-angled Artin group embeds into a right-angled Artin group whose defining graph is the complement of a finite interval graph. For example, according to \cite{MR3072113}, we already know that, 
\begin{itemize}
	\item for every $n \geq 1$, $A(C_n) \hookrightarrow A(C_5) \hookrightarrow A(\overline{P_7})$,
	\item for any finite tree $T$, $A(T) \hookrightarrow A(C_5) \hookrightarrow A(\overline{P_7})$,
\end{itemize}
where $P_7$ is an interval graph.

\medskip \noindent
Finally, it is worth noticing that the results proved in this paper depend on the chosen semigroup presentation $\mathcal{P}$ associated to the diagram group $D(\mathcal{P},w)$. Therefore, an interesting problem would be to find algebraic properties allowing us to choose $\mathcal{P}$ with suitable finiteness properties. For example:

\begin{question}
Let $G$ be a diagram group. When does there exist a semigroup presentation $\mathcal{P}= \langle \Sigma \mid \mathcal{R} \rangle$ and a base word $w \in \Sigma^+$ such that $G \simeq D(\mathcal{P},w)$ with $[w]_{\mathcal{P}}$ finite or $S(\mathcal{P},w)$ finite-dimensional?
\end{question}

\noindent
We suspect that the problem above could be solved by considering cohomological finiteness conditions. In particular, we do not know examples of finitely-presented diagram groups of finite algebraic dimension which cannot be expressed as $D(\mathcal{P},w)$ with $[w]_{\mathcal{P}}$ finite. (Without the hypothesis of being finitely-presented, $\mathbb{F}_{\infty}$ and $\mathbb{Z} \bullet \mathbb{Z}$ would be counterexamples, although the associated Squier complexes may be chosen of finite dimension.)

\addcontentsline{toc}{section}{References}

\bibliographystyle{alpha}
\bibliography{HSCC}

\def\polhk#1{\setbox0=\hbox{#1}{\ooalign{\hidewidth
  \lower1.5ex\hbox{`}\hidewidth\crcr\unhbox0}}}
\begin{thebibliography}{CRDK13}

\bibitem[Ago13]{MR3104553}
Ian Agol.
\newblock The virtual {H}aken conjecture.
\newblock {\em Doc. Math.}, 18:1045--1087, 2013.
\newblock With an appendix by Agol, Daniel Groves, and Jason Manning.

\bibitem[AGS06]{MR2271228}
G.~N. Arzhantseva, V.~S. Guba, and M.~V. Sapir.
\newblock Metrics on diagram groups and uniform embeddings in a {H}ilbert
  space.
\newblock {\em Comment. Math. Helv.}, 81(4):911--929, 2006.

\bibitem[Bau81]{MR634562}
A.~Baudisch.
\newblock Subgroups of semifree groups.
\newblock {\em Acta Math. Acad. Sci. Hungar.}, 38(1-4):19--28, 1981.

\bibitem[BDS07]{MR2308852}
Sergei Buyalo, Alexander Dranishnikov, and Viktor Schroeder.
\newblock Embedding of hyperbolic groups into products of binary trees.
\newblock {\em Invent. Math.}, 169(1):153--192, 2007.

\bibitem[BH99]{MR1744486}
Martin~R. Bridson and Andr{\'e} Haefliger.
\newblock {\em Metric spaces of non-positive curvature}, volume 319 of {\em
  Grundlehren der Mathematischen Wissenschaften [Fundamental Principles of
  Mathematical Sciences]}.
\newblock Springer-Verlag, Berlin, 1999.

\bibitem[B{\'S}99]{MR1703363}
Werner Ballmann and Jacek {\'S}wi{\polhk{a}}tkowski.
\newblock On groups acting on nonpositively curved cubical complexes.
\newblock {\em Enseign. Math. (2)}, 45(1-2):51--81, 1999.

\bibitem[CD95]{MR1303028}
Ruth Charney and Michael~W. Davis.
\newblock The {$K(\pi,1)$}-problem for hyperplane complements associated to
  infinite reflection groups.
\newblock {\em J. Amer. Math. Soc.}, 8(3):597--627, 1995.

\bibitem[CRDK13]{MR3072113}
Montserrat Casals-Ruiz, Andrew Duncan, and Ilya Kazachkov.
\newblock Embedddings between partially commutative groups: two
  counterexamples.
\newblock {\em J. Algebra}, 390:87--99, 2013.

\bibitem[CSS08]{MR2422070}
John Crisp, Michah Sageev, and Mark Sapir.
\newblock Surface subgroups of right-angled {A}rtin groups.
\newblock {\em Internat. J. Algebra Comput.}, 18(3):443--491, 2008.

\bibitem[Far03]{MR1978047}
Daniel~S. Farley.
\newblock Finiteness and {$\rm CAT(0)$} properties of diagram groups.
\newblock {\em Topology}, 42(5):1065--1082, 2003.

\bibitem[Gen]{arXiv:1505.02053}
A.~Genevois.
\newblock Hyperbolic diagram groups are free.
\newblock {\em arXiv:1505.02053}.

\bibitem[GS97]{MR1396957}
Victor Guba and Mark Sapir.
\newblock Diagram groups.
\newblock {\em Mem. Amer. Math. Soc.}, 130(620):viii+117, 1997.

\bibitem[GS99]{MR1725439}
V.~S. Guba and M.~V. Sapir.
\newblock On subgroups of the {R}. {T}hompson group {$F$} and other diagram
  groups.
\newblock {\em Mat. Sb.}, 190(8):3--60, 1999.

\bibitem[GS06a]{MR2193190}
V.~S. Guba and M.~V. Sapir.
\newblock Diagram groups and directed 2-complexes: homotopy and homology.
\newblock {\em J. Pure Appl. Algebra}, 205(1):1--47, 2006.

\bibitem[GS06b]{MR2193191}
V.~S. Guba and M.~V. Sapir.
\newblock Diagram groups are totally orderable.
\newblock {\em J. Pure Appl. Algebra}, 205(1):48--73, 2006.

\bibitem[Hag08]{MR2413337}
Fr{\'e}d{\'e}ric Haglund.
\newblock Finite index subgroups of graph products.
\newblock {\em Geom. Dedicata}, 135:167--209, 2008.

\bibitem[Hum]{arXiv:1207.2132}
D.~Hume.
\newblock Embedding mapping class groups into finite products of trees.
\newblock {\em arXiv:1207.2132}.

\bibitem[HW08]{MR2377497}
Fr{\'e}d{\'e}ric Haglund and Daniel~T. Wise.
\newblock Special cube complexes.
\newblock {\em Geom. Funct. Anal.}, 17(5):1551--1620, 2008.

\bibitem[HW15]{MR3320891}
Mark~F. Hagen and Daniel~T. Wise.
\newblock Cubulating hyperbolic free-by-cyclic groups: the general case.
\newblock {\em Geom. Funct. Anal.}, 25(1):134--179, 2015.

\bibitem[Kil97]{MR1448329}
Vesna Kilibarda.
\newblock On the algebra of semigroup diagrams.
\newblock {\em Internat. J. Algebra Comput.}, 7(3):313--338, 1997.

\bibitem[LW13]{MR3118410}
Joseph Lauer and Daniel~T. Wise.
\newblock Cubulating one-relator groups with torsion.
\newblock {\em Math. Proc. Cambridge Philos. Soc.}, 155(3):411--429, 2013.

\bibitem[MM99]{MR1672910}
Terry~A. McKee and F.~R. McMorris.
\newblock {\em Topics in intersection graph theory}.
\newblock SIAM Monographs on Discrete Mathematics and Applications. Society for
  Industrial and Applied Mathematics (SIAM), Philadelphia, PA, 1999.

\bibitem[MS13]{MR3073916}
John~M. Mackay and Alessandro Sisto.
\newblock Embedding relatively hyperbolic groups in products of trees.
\newblock {\em Algebr. Geom. Topol.}, 13(4):2261--2282, 2013.

\bibitem[NR98]{MR1459140}
Graham~A. Niblo and Martin~A. Roller.
\newblock Groups acting on cubes and {K}azhdan's property ({T}).
\newblock {\em Proc. Amer. Math. Soc.}, 126(3):693--699, 1998.

\bibitem[NR03]{MR1983376}
G.~A. Niblo and L.~D. Reeves.
\newblock Coxeter groups act on {${\rm CAT}(0)$} cube complexes.
\newblock {\em J. Group Theory}, 6(3):399--413, 2003.

\bibitem[Pau01]{MR1883722}
Scott~D. Pauls.
\newblock The large scale geometry of nilpotent {L}ie groups.
\newblock {\em Comm. Anal. Geom.}, 9(5):951--982, 2001.

\bibitem[Sag95]{MR1347406}
Michah Sageev.
\newblock Ends of group pairs and non-positively curved cube complexes.
\newblock {\em Proc. London Math. Soc. (3)}, 71(3):585--617, 1995.

\bibitem[Sag14]{SageevCAT(0)}
Michah Sageev.
\newblock Cat(0) cube complexes and groups.
\newblock In {\em Geometric Group Theory}, volume~21 of {\em IAS/Park City
  Math. Ser.}, pages 7--54. 2014.

\bibitem[Ser03]{MR1954121}
Jean-Pierre Serre.
\newblock {\em Trees}.
\newblock Springer Monographs in Mathematics. Springer-Verlag, Berlin, 2003.
\newblock Translated from the French original by John Stillwell, Corrected 2nd
  printing of the 1980 English translation.

\bibitem[SW79]{MR564422}
Peter Scott and Terry Wall.
\newblock Topological methods in group theory.
\newblock In {\em Homological group theory ({P}roc. {S}ympos., {D}urham,
  1977)}, volume~36 of {\em London Math. Soc. Lecture Note Ser.}, pages
  137--203. Cambridge Univ. Press, Cambridge-New York, 1979.

\bibitem[Wis04]{MR2053602}
D.~T. Wise.
\newblock Cubulating small cancellation groups.
\newblock {\em Geom. Funct. Anal.}, 14(1):150--214, 2004.

\end{thebibliography}

\end{document}